\documentstyle[amssymb,12pt]{amsart}
\newtheorem{prop}{Proposition}[section]

\newtheorem{lemma}{Lemma}[section]
\newtheorem{thm}{Theorem}[section]
\newtheorem{corollary}{Corollary}[section]

\theoremstyle{remark}
\newtheorem{remark}{Remark}
\newtheorem{notation}{Notation}

\begin{document}
\newcommand{\nc}{\newcommand} \nc{\on}{\operatorname}
\nc{\pa}{\partial}
\nc{\cA}{{\cal A}}\nc{\cB}{{\cal B}}\nc{\cC}{{\cal C}}
\nc{\cE}{{\cal E}}\nc{\cG}{{\cal G}}\nc{\cH}{{\cal H}}
\nc{\cX}{{\cal X}}\nc{\cR}{{\cal R}}\nc{\cL}{{\cal L}}
\nc{\cK}{{\cal K}}
\nc{\sh}{\on{sh}}\nc{\Id}{\on{Id}}\nc{\Diff}{\on{Diff}}
\nc{\ad}{\on{ad}}\nc{\Der}{\on{Der}}\nc{\End}{\on{End}}
\nc{\res}{\on{res}}\nc{\ddiv}{\on{div}}\nc{\card}{\on{card}}
\nc{\Jac}{\on{Jac}}
\nc{\Imm}{\on{Im}}\nc{\limm}{\on{lim}}\nc{\Ad}{\on{Ad}}
\nc{\ev}{\on{ev}}
\nc{\Hol}{\on{Hol}}\nc{\Det}{\on{Det}}
\nc{\de}{\delta}\nc{\si}{\sigma}\nc{\ve}{\varepsilon}
\nc{\al}{\alpha}
\nc{\CC}{{\mathbb C}}\nc{\ZZ}{{\mathbb Z}}
\nc{\NN}{{\mathbb N}}\nc{\zz}{{\mathbf z}}
\nc{\AAA}{{\mathbb A}}\nc{\cO}{{\cal O}} \nc{\cF}{{\cal F}}
\nc{\la}{{\lambda}}\nc{\G}{{\mathfrak g}}
\nc{\A}{{\mathfrak a}}
\nc{\HH}{{\mathfrak h}}
\nc{\N}{{\mathfrak n}}\nc{\B}{{\mathfrak b}}
\nc{\La}{\Lambda}
\nc{\g}{\gamma}\nc{\eps}{\epsilon}\nc{\wt}{\widetilde}
\nc{\wh}{\widehat}
\nc{\bn}{\begin{equation}}\nc{\en}{\end{equation}}
\nc{\SL}{{\mathfrak{sl}}}

%
%
%

\newcommand{\ldar}[1]{\begin{picture}(10,50)(-5,-25)
\put(0,25){\vector(0,-1){50}}
\put(5,0){\mbox{$#1$}} 
\end{picture}}

\newcommand{\lrar}[1]{\begin{picture}(50,10)(-25,-5)
\put(-25,0){\vector(1,0){50}}
\put(0,5){\makebox(0,0)[b]{\mbox{$#1$}}}
\end{picture}}

\newcommand{\luar}[1]{\begin{picture}(10,50)(-5,-25)
\put(0,-25){\vector(0,1){50}}
\put(5,0){\mbox{$#1$}}
\end{picture}}

\title[Some examples of quantum groups in higher genus]
{Some examples of quantum groups in higher genus}

\author{B. Enriquez}

\address{B.E.:Centre de Math\'ematiques, Ecole Polytechnique, 
URA 169 du CNRS, 91128 Palaiseau, France}

\address{FIM, ETH-Zentrum, HG G46, CH-8092 Z\"urich, Switzerland}

\author{V. Rubtsov}

\address{V.R.: D\'ept. de Math\'ematiques, Univ. d'Angers, 
2, Bd. Lavoisier, 49045 Angers, France}

\address{ITEP, 25 Bol. Cheremushkinskaia, 117259 Moscou, Russia}

\date{december 1997}

\begin{abstract}
  This is a survey of our construction of current algebras, associated
  with complex curves and rational differentials.  We also study in
  detail two classes of examples. The first is the case of a rational
  curve with differentials $z^n dz$; these algebras are ``building
  blocks'' for the quantum current algebras introduced in our earlier
  work. The second is the case of a genus $>1$ curve $X$, endowed with
  a regular differential having only double zeroes.
\end{abstract}

\maketitle

In our papers \cite{Enr-Rub,ER:qH}, we introduced a family of 
quasi-Hopf algebras, associated with complex curves and rational 
differentials. These algebras are the quantizations of quasi-Lie 
bialgebra structures that had been defined by V. Drinfeld in 
\cite{Dr}. 



Our purpose here is to first present (sect. \ref{general}) 
a survey of the constructions of
\cite{Enr-Rub,ER:qH}. After this, we present some examples.
 
Let us first recall the construction 
of the quasi-Lie bialgebras of \cite{Dr}. We fix a semisimple
Lie algebra $\bar\G$, a rational 
curve $X$ and a nonzero rational form $\omega$ on it. We denote
by $S$ a set of points of $X$, containing all poles and zeroes 
of $\omega$, and by $\cK_S$ the direct sum of local fields of $X$
at the points of $S$. We define $R$ as the subring of $\cK_S$ 
formed by the Laurent expansions at the points of $S$ of the 
regular functions on $X-S$. Then $\cK_S$ is endowed with a 
scalar product defined by $\omega$, and $R$ is then a maximal 
isotropic subspace of $\cK_S$. We fix a maximal isotropic 
supplementary $\La$ to $R$ in $\cK_S$. The Lie algebra 
$\bar\G \otimes \cK_S$ being endowed with a product pairing, 
we then have a direct sum decomposition 
$$
\bar\G\otimes\cK_S = (\bar\G\otimes R) 
\oplus (\bar\G\otimes\La) 
$$
in isotropic subspaces, whose first summand is a Lie subalgebra. 
This defines quasi-Lie bialgebra structures on 
$\bar\G\otimes\cK_S$ and $\bar\G\otimes R$: we have maps
$\delta_{\cK_S}$ and $\delta_{R}$ from $\bar\G\otimes\cK_S$ 
and $\bar\G\otimes R$ to their second exterior
power, and an element $\phi$ of $\wedge^3(\bar\G\otimes R)$, 
satisfying some compatibility contitions (see \cite{Dr}). 
The problem of their quantization is to construct quasi-Hopf 
algebras $U_{\hbar}(\bar\G\otimes \cK_S)$  and 
$U_{\hbar}(\bar\G\otimes R)$ deforming their enveloping
algebras, coproducts on these algebras deforming 
extensions of $\delta_{\cK_S}$ and of $\delta_{R}$, and 
elements of the tensor cubes of these algebras, satisfying the 
quasi-Hopf algebra axioms.


In \cite{Enr-Rub}, we solved this problem for double 
extensions of these quasi-Lie bialgebra structures (these
extensions are defined by the usual cocycle on current algebras, 
and by a derivation), in the special case when $\bar\G= \SL_2$. 
For this, we defined semi-infinite twists
of these structures, in the spirit of Drinfeld's new realizations, 
and quantized them (see \cite{new-real,Enr-Rub}). These 
quantizations, $(U_{\hbar}\G,\Delta)$ and $(U_{\hbar}\G,\bar\Delta)$, 
are related by some twist operation. 
In \cite{ER:qH}, we constructed in $U_{\hbar}\G$, a 
subalgebra $U_{\hbar}\G_R$ deforming $U\G_R$. Using the coideal 
properties of this algebra with respect to $\Delta$ and 
$\bar\Delta$, we reduced the problem of finding a quasi-Hopf
structure on $U_{\hbar}\G$, preserving $U_{\hbar}\G_R$, to some
decomposition problem on the twist $F$; this problem was solved
using some results on Hopf duality.  We close the section by a result
(Prop. \ref{simply:laced} and Cor. \ref{van:q(z,w)}) characterizing
the zeroes and poles of the structure function $q(z,w)$ defining
$U_{\hbar}\G$.  

In section \ref{level:0}, we come back on the quantization problem 
for the non-extended quasi-Lie bialgebra structures. We remark 
that there, due to the absence of derivation, much more relations
are possible for the algebra $U_{\hbar}\G$. We describe these 
relations, and how the construction of quasi-Hopf algebras described
above can be generalized in that situation. 



Let us say here some words on the applications of the 
constructions of sect. \ref{general} in genera $0$ and $1$. 
In genus zero, and 
with $\omega = dz$, this construction agrees with the quantum 
currents presentation of double Yangians; in 
\cite{EF:rat}, we derived another expression of Khoroshkin-Tolstoy 
twists (\cite{KhT}) relating Drinfeld's coproduct for the double
Yangians with the usual ($L$-operator) one. 


In the elliptic case, with $\omega$ regular, these algebras
were related with Felder's elliptic quantum groups (\cite{ell:QG}). 
In both situations, the problem was to suitably refine the 
decomposition of $F$, using additional conditions provided by 
algebras ``opposite'' the the regular one. 

Our goal in the next sections is to present some other examples
of the construction of sects. \ref{general}, \ref{level:0}. 

In sect. \ref{rat:sect}, we treat the case of a rational curve
with differential $z^n dz$. The corresponding algebra is denoted
$U_{z^{n}dz}\G$. We compute its structure
coefficients $q(z,w)$ in Prop. \ref{nesher}. We also study 
some properties of $U_{z^{n}dz}\G$: making use of the 
$\ZZ_{n+1}$-symmetry of the situation, we give a presentation 
of the vertex relations of $U_{z^{n}dz}\G$ not involving the 
$n+1$-st roots as in 
the expression of $q(z,w)$, which allows to make sense of this
algebra of vertex relations for complex values of the deformation 
parameter, 
without any completion procedure. We then show that, like what 
happens for the Yangian algebra, the quantum algebras of this
family are isomorphic for all nonzero values of the quantum 
parameter.  

In Thm. \ref{isomorphisms}, we show that the algebras 
$U_{z^n dz}\G$ are ``building blocks'' for the quantum current 
algebras $U_{\hbar}\G$ constructed in sect. \ref{general} -- that
is, each algebra $U_{\hbar}\G$ is isomorphic to a tensor product
of algebras  $U_{z^n dz}\G$ with their centers identified.  

In sect. \ref{g>1}, we turn to the case of a curve $X$ of genus $>1$, 
with a differential $\omega$ regular and having 
only double zeroes. The existence of such a form is a well-known 
fact from Riemann surface theory (see e.g. \cite{Mum,Fay}), and 
is equivalent to the existence of odd theta-characteristics.  

In this situation, we first construct isotropic supplementaries in 
local fields, 
using functions on $X$, which are multivalued along the $b$-cycles
(sect. \ref{suppl:sect}). We then compute the Green kernel of this
decomposition (sect. \ref{green:sect}). Applying results of sect. 
\ref{level:0}, we then present relations for quantizations of the 
centerless versions of these algebras (sect. 
\ref{level:0:theta:rels}). 
We remark that after a finite twist, the above decomposition  has a 
part consisting of the regular functions on $X$ minus some points. 
We derive from this the construction of a regular subalgebra in the
quantum currents algebra (sect. \ref{regular:level0:sect}). We show
(Rem. \ref{g>1:nonzero:level}) how these results may be extended for 
doubly extended quasi-Lie bialgebra structures. 

The zero-level relations might have special interest for the shift 
parameter $\hbar h$ (which belongs to the Jacobian of $X$) nonformal 
and torsion, like what happens for elliptic quantum 
groups with torsion parameter (\cite{Fel-Varch}). 

In both algebras, we also construct deformations of the 
enveloping algebra of $\bar\G \otimes (\cK_\delta \oplus \cO_{S'})$  
-- here $\delta$ is the set of zeroes of $\omega$, $S'$ are 
other marked points on $X$, and $\cK_\delta$ and $\cO_{S'}$ are the
sums of local fields and rings at these points. This might be useful 
for constructing induced modules. 

Let us now discuss possible prolongations of the present work. 
We did not examine degenerations of our constructions; 
of special interest should be rational or elliptic curves 
with points identified. In particular, in the latter
case, and with $\omega = dz$, the two types of algebras studied 
here (which are defined either by applying $q^\pa$ to variables 
in the structure relations, or by some shifts in theta-functions) 
should coincide. 

We also hope that Cor. \ref{van:q(z,w)} will help 
to treat quantum conformal blocks in higher genus as it
was done by B. Feigin and A. Stoyanovsky in \cite{FS} in the case
of the affine Kac-Moody algebra $\hat\SL_2$, 
and in \cite{Feigin-Ding} in the quantum affine case. 

Finally, we also would like to mention the papers 
\cite{FO,Feigin-Miwa}, which should be closely connected 
with some problems left open in \cite{ER:qH}: 
identification of the twists of twists of \cite{EF:rat} 
with those of Khoroshkin-Tolstoy in \cite{KhT}; and 
construction of quantum Serre relations in higher genus.  

We would like to thank J. Ding, B. Feigin, G. Felder, K. Gawedzki, 
K. Hasegawa, M. Jimbo, Y. Kosmann-Schwarzbach, H. Konno and 
J. Shiraishi for discussions 
related to the subject of this paper. B.E. would like to thank 
T. Miwa for 
invitation to RIMS, where some part of this work was carried out 
in a very stimulating atmosphere. V.R. would also like to
thank M. Audin and M. Rosso for their invitation to ULP Strasbourg-1 
university; he also acknowledges support of grants RFFI-97-01-01101 and 
INTAS-96-196.

\section{Review of quantum current and quasi-Hopf 
algebras in higher genus} 
\label{general}

\subsection{Manin pairs, triples and classical twists} \label{def:sect}

Let $X$ be a smooth, connected, compact complex curve, and $\omega$ be a 
nonzero
meromorphic differential on $X$. Let $S$ be a finite set of points of
$X$, containing the set $S_{0}$ of its zeros and
poles. For each $s\in S$, let $\cK_{s}$ be the local field at $s$ and
$$
\cK=\oplus_{s\in S}\cK_{s}. 
$$
Let $R$ be the ring of meromorphic functions on $X$, regular outside
$S$; $R$ can be viewed as a subring of $\cK$. $R$ is endowed with the
discrete topology and $\cK$ with its usual (formal series) topology. 
Let us define on $\cK$ the bilinear form 
$$
\langle f,g \rangle_{\cK}=\sum_{s\in S}\res_{s}(fg\omega), 
$$
and the derivation 
$$
\pa f= df/\omega.
$$
We will use the notation ${\mathfrak x}(A)={\mathfrak x}\otimes A$, for any ring
$A$ over $\CC$ and complex Lie algebra $\mathfrak x$. 

Let $\bar\G=\SL_2(\CC)$. Define on $\bar\G(\cK)$ the bilinear form
$\langle, \rangle_{\bar\G(\cK)}$ by 
$$ 
\langle x\otimes \eps, y\otimes \eta\rangle_{\bar\G(\cK)}=\langle
x,y\rangle_{\bar\G}\langle \eps, \eta \rangle_{\cK}
$$
for $x,y\in\bar\G, \eps,\eta\in \cK$, 
$\langle, \rangle_{\bar\G}$ being the
Killing form of $\bar\G$, the derivation $\pa_{\bar\G(\cK)}$ 
by $\pa_{\bar\G(\cK)}(x\otimes \eps)=x\otimes
\pa\eps$, for $x\in \bar\G, \eps\in \cK$, 
and the cocycle
$$
c(\xi,\eta)=\langle \xi,\pa_{\bar\G(\cK)} \eta\rangle_{\bar\G(\cK)}. 
$$
Let $\hat{\G}$ be the central extension of 
$\bar\G(\cK)$ by this cocycle. We
then have 
$$
\hat{\G}=\bar\G(\cK)\oplus\CC K, 
$$
with bracket such that $K$ is central, and  $[\xi,\eta]
=([\bar \xi,\bar \eta],c(\bar \xi,\bar \eta)K)$, for any $\xi,\eta\in 
\hat{\G}$ with
first components $\bar \xi,\bar \eta$. 

Let us denote by $\pa_{\hat{\G}}$ the derivation of 
$\hat{\G}$ defined by $\pa_{\hat{\G}}(\xi,0)=(\pa_{\bar\G(\cK)}\xi,0)$ and
$\pa_{\hat{\G}}(K)=0$. 

Let $\G$ be the skew product of $\hat{\G}$ with $\pa_{\hat{\G}}$. We
have 
$$
\G=\hat{\G}\oplus \CC D, 
$$
with bracket such that $\hat{\G}\to \G$, $\xi\mapsto (\xi,0)$ is a Lie
algebra morphism, and $[D,(\xi,0)]=(\pa_{\hat{\G}}(\xi),0)$ for
$\xi\in\hat{\G}$. 

View $\bar\G(\cK)$ as a subspace of 
$\G=\hat{\G}\oplus \CC D=\bar\G(\cK)\oplus \CC K
\oplus \CC D$, by $\xi\mapsto (\xi,0,0)$. 
Define on $\G$ the pairing $\langle, \rangle_{\G}$ by $\langle
K,D\rangle_{\G} =1$, $\langle K, \bar\G(\cK)\rangle_{\G} = \langle D,
\bar\G(\cK)\rangle_{\G} 
=0$, $\langle \xi,\eta \rangle_{\G}=\langle \xi,\eta 
\rangle_{\bar\G(\cK)}$
for $\xi,\eta\in \bar\G(\cK)$. 

Endow $\bar\G(\cK)$ with $\langle, \rangle_{\bar\G(\cK)}$. The subspace
$\bar\G(R)\subset \bar\G(\cK)$
is a maximal isotropic subalgebra of $\bar\G(\cK)$. 
Drinfeld's Manin pair is $(\bar\G(\cK), \bar\G(R))$ (see \cite{Dr}). 
In \cite{Enr-Rub}, we introduced the following extension of this pair. 
Let $\G_{R}=\bar\G(R)\oplus\CC D$; $\G_{R}\subset \G$ is a maximal isotropic
subalgebra of $\G$. The extended Drinfeld's Manin pair of 
\cite{Enr-Rub} is then $(\G,\G_{R})$. 

In \cite{Enr-Rub}, we also introduced the following
Manin triple. Let $\La$ be a Lagrangian complement to $R$ in $\cK$,
commensurable with $\oplus_{s\in S}\cO_{s}$ (where $\cO_{s}$ is the
completed local ring at $s$). 
Let $\N_{+}=\CC
e$, $\N_{-}=\CC f$, $\HH=\CC h$. Let 
$$
\G_{+}=\HH(R)\oplus \N_{+}(\cK)\oplus \CC D, 
\quad \G_{-}=(\HH\otimes \La)\oplus \N_{-}(\cK) \oplus \CC K,
$$
then $\G=\G_{+}\oplus \G_{-}$, and 
both $\G_{+}$ and $\G_{-}$ are maximal isotropic subalgebras of 
$\G$. The Manin triple is then $(\G,\G_{+},\G_{-})$. 

We will also consider the following Manin triple, that we may consider
as being obtained from the previous one by the action of the nontrivial
element of the Weyl group of $\bar\G$. 
Let 
$$
\bar{\G}_{+}=\HH(R)\oplus \N_{-}(\cK)\oplus \CC D, 
\quad \bar{\G}_{-}=(\HH\otimes \La)\oplus \N_{+}(\cK) \oplus \CC K,
$$
then $(\G,\bar{\G}_{+},\bar{\G}_{-})$ again forms a Manin triple. 
\medskip

According to \cite{QG}, to each of the Manin triples
$(\G,\G_{+},\G_{-})$ and $(\G,\bar{\G}_{+},\bar{\G}_{-})$ is
associated a Lie bialgebra structure on $\G$; denote by
$\delta,\bar\delta:\G\to \G\hat\otimes\G$ the corresponding cocycle maps.

Let $\G_{\La} = (\bar\G\otimes \La)\oplus \CC K \subset \G$; $\G_{\La}$ is a
Lagrangian complement of $\G_{R}$ in $\G$. It induces a Lie
quasi-bialgebra structure on $\G_{R}$, and 
from \cite{Kosmy} follows also that there is a Lie quasi-bialgebra
structure on $\G$, associated to the Manin pair $(\G,\G_{R})$ and to
$\G_{\La}$; we denote
by $\delta_{R}:\G\to\G\hat\otimes\G$ the corresponding cocycle map. 

These Lie (quasi-)bialgebra structures on $\G$ are related by the
following classical twist operations. 

Let $(e^{i})_{i\in \NN}, (e_{i})_{i\in \NN}$ be dual bases of $R$ and
$\La$; we choose them is such a way that $e_{i}$ tends to $0$ when $i$
tends to $\infty$.
Let $\eps^{i}, \eps_{i}, i\in \ZZ$ be dual bases of $\cK$, defined by
$\eps_{i}=e_{i}, \eps^{i}=e^{i}, i\ge 0$,  $\eps_{i}=e^{-i-1},
\eps^{i}=e_{-i-1}, i< 0$.

\begin{lemma} (see \cite{ER:qH})
Let $f=\sum_{i\in\ZZ}e[\eps^{i}]\otimes f[\eps_{i}]$; 
$f=f_{1}+f_{2}$, with 
$$
f_{1}=\sum_{i\in\NN}e[e_{i}]\otimes f[e^{i}], 
$$ 
and 
$$
f_{2}=\sum_{i\in\NN}e[e^{i}]\otimes f[e_{i}]. 
$$ 
For $\xi\in\G$, we have
$$
\delta_{R}(\xi)=\delta(\xi)+[f_{1}, \xi\otimes 1 + 1\otimes \xi], 
\quad 
\bar\delta(\xi)=\delta_{R}(\xi)+[f_{2}, \xi\otimes 1 + 1\otimes \xi]. 
$$
\end{lemma}

\subsection{Results on Green kernels}

\begin{notation}
For $a= a(z,w)$ a function of two variables $z,w$, 
we denote by $a^{(21)}$ the function $a^{(21)}(z,w) = a(w,z)$. 
\end{notation}

Let us fix dual bases  
$(e^{i})_{i\in \NN}, (e_{i})_{i\in \NN}$ of $R$ and
$\La$. 
Let $G\in R\hat\otimes \La$ be the series 
$$
G=\sum_{i}e^{i}\otimes e_{i}; 
$$
it is called the Green kernel of $(X,\omega,S,\La)$
Note that $R\hat\otimes k$ is an algebra, to which $G$ belongs.
Let 
$$
\gamma=(\pa\otimes 1)G-G^{2}; 
$$
then 
\begin{lemma} \label{barbara} (see \cite{Enr-Rub})
$\gamma$ belongs to $R\otimes R$. 
\end{lemma}

Let $\hbar$ be a formal variable and let $T:k[[\hbar]]\to k[[\hbar]]$
be the operator equal to 
$$
T={{\sh(\hbar \partial)} \over {\hbar \partial}}.
$$
We will use the notation $q = e^{\hbar}$. 

Let $(\gamma_{i})_{i\ge 0}$ be a set of free variables, and
$\phi,\psi\in \hbar\CC[\gamma_{i}] [[\hbar]]$ be the solutions of
$$
{{\pa\psi}\over{\pa\hbar}} 
= D \psi -1-\gamma_{0}\psi^{2}, \quad
{{\pa\phi}\over{\pa\hbar}} 
=D \phi - \gamma_{0} \psi , 
$$
where $D = \sum_{i\ge 0}\gamma_{i+1}{\pa\over{\pa \gamma_{i}}}$; 
then 
\begin{prop} (see \cite{Enr-Rub}, Prop. 3)
$$
\sum_{i\in\NN}{{q^{\pa} - 1}\over{\pa}}e^{i} \otimes e_{i}
=
\phi(\hbar , (\pa^{i}\otimes 1)\gamma)
- \ln (1 + G \psi(\hbar , (\pa^{i}\otimes 1)\gamma)). 
$$
\end{prop}

{}From this Prop. follows: 

\begin{prop} (see \cite{Enr-Rub})
For certain elements $\phi\in (R\otimes R)[[\hbar]]$, 
$\psi_{+},\psi_{-}\in \hbar (R\otimes R)[[\hbar]]$, we
have the following identities in $(R\hat\otimes \cK)[[\hbar]]$
$$
\sum_{i} Te^{i}\otimes e_{i}= \phi+{1\over{2\hbar}}\ln{{1+G\psi_{-}}
\over{1+G\psi_{+}}},
\quad \sum_{i} e^{i}\otimes Te_{i}= 
-{\phi}^{(21)}+{1\over{2\hbar}}\ln{{1-G{\psi}^{(21)}_{+}}
\over{1-G{\psi}^{(21)}_{-}}}
.  
$$
\end{prop}

\begin{lemma} (see \cite{Enr-Rub}) The expression 
$\sum_{i}Te^{i}\otimes e_{i} - e^{i}\otimes Te_{i}$
belongs to $S^{2}(R)[[\hbar]]$. We will denote by $\tau$ any element
of $(R\otimes R)[[\hbar]]$, such that 
\begin{equation} \label{id-tau}
\tau+{\tau}^{(21)}=\sum_{i}Te^{i}\otimes e_{i} - e^{i}\otimes 
Te_{i}. 
\end{equation}
\end{lemma}

Note that $\sum_{i}Te^{i}\otimes e_{i}$ is well-defined in
$(R\hat\otimes \cK)[[\hbar]]$, because 
$e_{i}$ tends to zero as $i$ tends
to infinity. Since $\pa$ is a continuous map from $\cK$ to itself, the
same is true for the sequence $\pa^{k}e_{i}$. So $\sum_{i}e^{i}\otimes
Te_{i}$ is well-defined in the same space; 
$\sum_{i\in \ZZ}T\eps^{i}\otimes \eps_{i} -\eps^{i}\otimes
T\eps_{i}$
is well-defined in 
$(\cK\bar\otimes \cK)[[\hbar]]$ for the same reasons.

Define 
\begin{equation} \label{q(z,w)}
q(z,w) = q^{2(\tau-\phi)} {{1 + \psi_+ G}\over{1+\psi_- G}}(z,w).  
\end{equation}

\subsection{Hopf algebras $(U_{\hbar}\G,\Delta)$,
$(U_{\hbar}\G,\bar\Delta)$ and the twist connecting them}

In \cite{Enr-Rub}, we introduced a Hopf algebra 
$(U_{\hbar}\G,\Delta)$ quantizing $(\G,\delta)$. 

It is the quotient of $T(\G)\hat{}[[\hbar]]$ by the following
relations. Let $e,f,h$ be the Chevalley basis of $\SL_{2}(\CC)$.  
Denote in $T(\G)\hat{}[[\hbar]]$, the element $x\otimes \eps\in 
\bar\G(\cK)\subset
\G$ of
$\G$ by $x[\eps]$ and let for $r\in R$, $h^{+}[r]=h[r]$,
$h^{-}[\la]=h[\la]$. 
Introduce the generating series 
$$
e(z)=\sum_{i\in\ZZ}e[\eps_{i}]\eps^{i}(z), \quad
f(z)=\sum_{i\in \ZZ}f[\eps_{i}]\eps^{i}(z), 
$$
$$
h^{+}(z)=\sum_{i\in\NN}h^{+}[e^{i}]e_{i}(z), \quad
h^{-}(z)=\sum_{i\in\NN}h^{-}[e_{i}]e^{i}(z). 
$$
The Cartan fields are arranged in the series
$$
K^{+}(z) = e^{((T+U)h^+)(z)} , \quad K^{-}(z) = e^{-\hbar h^-(z)};   
$$
here $U$ is the linear operator from $\La$ to $R[[\hbar]]$ defined by 
$U(\la) = \langle \tau,1\otimes\la \rangle$.  
The relations for $U_{\hbar}\G$ are the coefficients of 
\begin{equation} \label{K-K}
[K^{+}(z),K^{+}(w)]=0, \quad
(K^{+}(z),K^{-}(w))
={{q(z,w)}\over{q(z,q^{-K\pa}(w))}},
\end{equation}
\begin{equation} \label{K-K-}
(K^{-}(z),K^{-}(w))={{q(q^{-K\pa}(z),q^{-K\pa}(w))}\over{q(z,w)}},
\end{equation}
\begin{equation}\label{K-e}
(K^{+}(z),e(w))= q(z,w), \quad (K^{-}(z),e(w))= q(w, q^{-K\pa}(z)), 
\end{equation}
\begin{equation}\label{K-f}
(K^{+}(z),f(w))= q(w,z), \quad (K^{-}(z),f(w))= q(z,w), 
\end{equation}
\begin{equation} \label{vx-e}
(z_s - w_s)(1+\psi_+ G)(z,w) e(z)e(w) = 
(z_s - w_s)e^{2(\tau-\phi)(z,w)}(1+\psi_- G)(z,w) e(w)e(z), 
\end{equation}
\begin{equation} \label{vx-f}
(z_s - w_s)e^{2(\tau-\phi)(z,w)}(1+\psi_- G)(z,w) f(z)f(w), 
= 
(z_s - w_s)(1+\psi_+ G)(z,w) f(w)f(z) 
\end{equation}
\begin{equation} \label{dsc-e-f}
[e(z),f(w)] = \delta(z,w)K^{+}(z) -\delta(z,q^{-K\pa}(w))
K^{-}(w)^{-1}; 
\end{equation}
\begin{equation}
[D,x^{\pm}(z)] = - (\pa x^{\pm})(z) + \hbar (A h^+)(z) x^{\pm}(z), 
\end{equation}
\begin{equation} \label{D-others}
[D,K^{\pm}(z)] = - (\pa K^{\pm})(z) + \hbar (B^{\pm} h^+)(z) 
K^{\pm}(z), 
\end{equation}
$K$ is central.  
We used the standard notation $(a,b)$ for the group commutator
$aba^{-1}b^{-1}$; we also set $\delta(z,w)= G(z,w) + G^{(21)}(z,w)$. 
Here $A$ and $B^{\pm}$ are operators from $\La$ to $R[[\hbar]]$; 
$A$ is defined by $A(\la) = T((\pa\la)_R) + \pa(U\la) - 
U((\pa\la)_{\La})$; formulas for $B^{\pm}$ can be extracted
from \cite{Enr-Rub}.

The formulas 
\begin{equation}  \label{Delta-K}
\Delta(K)=K\otimes 1 + 1 \otimes K
\end{equation}
\begin{equation}    \label{Delta-h}
\Delta(h^{+}[r])=h^{+}[r]\otimes 1+1\otimes h^{+}[r], \quad
\Delta(h^{-}(z))=h^{-}(z)\otimes 1+1\otimes (q^{-K_{1}\pa}h^{-})(z), 
\end{equation}
\begin{equation}   \label{Delta-e}
\Delta(e(z))=e(z)\otimes K^+(z) + 1\otimes e(z), 
\end{equation}
\begin{equation}   \label{Delta-f}
\Delta(f(z))=f(z)\otimes 1+ K^-(z)^{-1}\otimes (q^{-K_{1}\pa}f)(z),   
\end{equation}
\begin{equation} \label{Delta-D}
\Delta(D)= D \otimes 1 + 1 \otimes D + \sum_{i\in \NN}{\hbar \over 4}
h^{+}[e^{i}] \otimes h^{+}[A e_{i}],
\end{equation} 
$r\in R$, for the coproduct, 
\begin{equation} \label{counit}
\varepsilon(h^{+}[r])=\varepsilon(h^{-}[\la])=\varepsilon(x[\eps])=
\varepsilon(D)=\varepsilon(K)=0,
\end{equation}
$x=e,f$, $r\in R, \la\in\La, \eps\in \cK$, 
for the counit,
define a topological (with respect to the completion introduced above)
Hopf algebra structure on $U_{\hbar}\G$. 

The coalgebra structure of $U_{\hbar}\bar{\G}$  is defined by the
coproduct 
\begin{equation} \label{Delta'-K}
\bar\Delta(K)=K\otimes 1 + 1 \otimes K, 
\end{equation}
\begin{equation}    \label{Delta'-h}
\bar\Delta(h^{+}[r])=h^{+}[r]\otimes 1+1\otimes h^{+}[r], \quad
\bar\Delta(\bar h^{-}(z))=h^{-}(z)\otimes
1+1\otimes (\bar q^{-K_1\pa}h^{-})(z), 
\end{equation}
\begin{equation}   \label{Delta'-e}
\bar\Delta(e(z))=(q^{\bar K_{2}\pa}(e\otimes K^-))(z)
+ 1\otimes e(z), 
\end{equation}
\begin{equation}   \label{Delta'-f}
\bar\Delta(f(z))=f(z)\otimes 1+ K^+(z) \otimes f(z),
\end{equation}
\begin{equation} \label{Delta'-D}
\bar\Delta(\bar D)= \bar D \otimes 1 + 1 \otimes \bar D 
+ \sum_{i\in \NN}{\hbar \over 4}
h^{+}[e^{i}] \otimes h^{+}[A e_{i}],
\end{equation} 
$r\in R$, the counit
\begin{equation} \label{counit'}
\bar\varepsilon(h^{+}[r])=\bar \varepsilon(h^{-}[\la])=\bar
\varepsilon(x[\eps])=
\bar \varepsilon(D)=\bar \varepsilon(K)=0,
\end{equation}
$x=e,f$, $r\in R, \la\in\La, \eps\in \cK$,

\begin{thm} \label{def:quantum:algebra} The pairs  
$(U_\hbar\G,\Delta)$ and $(U_\hbar\G,\bar\Delta)$ defined by the
above relations are Hopf algebras
quantizing the Lie bialgebras $(\G,\delta)$ and $(\G,\bar\delta)$. 
\end{thm}

This result means in particular
that $U_{\hbar}\G$ is a flat deformation of the enveloping
algebra $U\G$. 
This follows from the following 
Poincar\'e-Birkhoff-Witt-type (PBW) result: 

\begin{prop}  \label{PBW} (see \cite{ER:qH}, Prop. 4.1). 
Let $A$ be an algebra with generators $x_n,x\in\ZZ$, 
generating series $x(z) = \sum_{n\in\ZZ} x_n z^{-n}$, and
relations  defined by the modes of
$$
(z-w + \sum_{i\ge 1}\hbar^i a_i(z,w)) x(z)x(w)
= 
(z-w + \sum_{i\ge 1}\hbar^i b_i(z,w)) x(w)x(z), 
$$
for $a_i$ and $b_i$ series of $\CC[[z,w]][z^{-1},w^{-1}]$. 
Set $a(z,w) = z-w + \sum_{i\ge 1}\hbar^i a_i(z,w)$, 
$b(z,w) = z-w + \sum_{i\ge 1}\hbar^i b_i(z,w)$. 
Then if the series $a_i$ and $b_i$ satisfy 
$a(z,w)a(w,z) = b(z,w)b(w,z)$, 
$A$ is a 
flat deformation of the symmetric algebra in the variables
$x_n,n\in\ZZ$.  
\end{prop}

Thm. \ref{def:quantum:algebra} follows from a double construction.
More precisely, define $U_{\hbar}\HH_R$ and $U_{\hbar}\HH_\La$ as the
subalgebras of $U_{\hbar}\G$ generated by $D$ and the $h^+[r], r\in
R$, resp. $K$ and the $h[\la],\la\in\La$; and define $U_\hbar\N_+$ and
$U_{\hbar}\N_-$ as the subalgebras of $U_\hbar\G$ generated by the
$e[\eps],\eps\in\cK$, resp. the $f[\eps],\eps\in\cK$. Set 
$$
U_\hbar\G_+ = U_\hbar\HH_R U_\hbar\N_+, \quad
U_\hbar\G_- = U_\hbar\HH_\La U_\hbar\N_-,
$$
and
$$
U_\hbar\bar\G_+ = U_\hbar\HH_R U_\hbar\N_-, \quad
U_\hbar\bar\G_- = U_\hbar\HH_\La U_\hbar\N_+. 
$$
Then: 
\begin{prop} (see \cite{Enr-Rub}) $U_\hbar\G_\pm$ and $U_\hbar\bar\G_\pm$
  are subalgebras of $U_\hbar\G$. $(U_\hbar\G_\pm,\Delta)$ and
  $(U_\hbar\bar\G_\pm,\bar\Delta)$ are Hopf subalgebras of
  $(U_\hbar\G,\Delta)$ and $(U_\hbar\G,\bar\Delta)$. Moreover,
  $(U_\hbar\G_+,\Delta)$ and $(U_\hbar\G_-,\Delta')$ are Hopf dual, as
  well as $(U_\hbar\bar\G_+,\bar\Delta)$ and
  $(U_\hbar\bar\G_-,\bar\Delta')$ and $(U_\hbar\G,\Delta)$ and
  $(U_\hbar\G,\bar\Delta)$ are the corresponding Drinfeld doubles.
\end{prop}

We then have: 

\begin{lemma} \label{id:pairing}
  The restriction to $U_\hbar\N_+\times U_\hbar\N_-$ of the Hopf
  pairing between $(U_\hbar\G_+,\Delta)$ and $(U_\hbar\G_-,\Delta')$
  agrees up to permutation of factors with the restriction to
  $U_\hbar\N_-\times U_\hbar\N_+$ of the Hopf pairing between
  $(U_\hbar\bar\G_+,\bar\Delta)$ and $(U_\hbar\bar\G_-,\bar\Delta')$.
\end{lemma}

Let us define the completion $U_{\hbar}\G \bar\otimes U_{\hbar}\G$
as follows. Let $I_{N}\subset U_{\hbar}\G$ be the left ideal generated
by the $x[\eps]$, $\eps\in \prod_{s\in
S}z_{s}^{N}\CC[[z_{s}]]$.  Define $U_{\hbar}\G \bar\otimes U_{\hbar}\G$
as the inverse limit of the $U_{\hbar}\G^{\otimes 2} / I_{N} \otimes 
U_{\hbar} \G + U_{\hbar}\G \otimes I_{N}$ (where the tensor products are
$\hbar$-adically completed). $U_{\hbar}\G \bar\otimes U_{\hbar}\G$ is
clearly a completion of $U_{\hbar}\G^{\hat\otimes 2}$. 

\begin{prop} Let $(\al^i),(\al_i)$ be bases of $U_\hbar\N_+$ and 
$U_\hbar\N_-$, dual for the pairing of Lemma \ref{id:pairing}. Set 
\begin{equation} \label{F}
F=\sum_i \al^i \otimes \al_i; 
\end{equation} 
$F$ belongs to $U_{\hbar}\G\bar\otimes U_{\hbar}\G$, and is a twist
tranforming $U_{\hbar}\G$ into $U_{\hbar}\bar{\G}$. More precisely,
\begin{equation} \label{twist}
\Ad (F)(\Delta(x))=\bar\Delta(x), 
\end{equation}
for each $x\in U_{\hbar} \G$. 
\end{prop}

\subsection{Universal $R$-matrices and Hopf algebra pairings}

Let $U_{\hbar}\G_{\pm}$ be the subalgebras of
$U_{\hbar}\G$ generated by $\G_{\pm}$. These are Hopf subalgebras of
$U_{\hbar}\bar{\G}$, dual to each other if $U_{\hbar}\G_{-}$ is endowed
with the coproduct opposite to $\Delta$.  
Then $U_{\hbar}\G_{+}$ and $U_{\hbar}\G_{-}$ are dual
to each other, and $U_{\hbar} \G$ is the corresponding double algebra. 

The pairing $\langle, \rangle_{U_{\hbar}\G}$ 
between $U_{\hbar}\G_{+}$ and $U_{\hbar}\G_{-}$ is defined
by 
$$
\langle h^{+}[r], h^{-}[\la]\rangle_{U_{\hbar}\G} 
={2\over \hbar}\langle r, \la\rangle_{k}, \quad
\langle e[\eps], f[\eta]\rangle_{U_{\hbar}\G} 
={1\over \hbar}\langle \eps, \eta\rangle_{k}, 
$$
for $\eps,\eta\in k$, $r\in R,\la\in\La$,
$$
\langle D,K\rangle_{U_{\hbar}\G} = 1, \quad 
\langle D, \bar\G(\cK) \rangle_{U_{\hbar}\G}=\langle K,
\bar\G(\cK)\rangle_{U_{\hbar}\G}=0, 
$$
and to be a Hopf algebra pairing, $U_{\hbar}\G_{-}$ being endowed the
coproduct opposite to the one given by its embedding in $U_{\hbar}\G$.

Let $U_{\hbar}\bar{\G}_{\pm}$ be the subalgebras of $U_{\hbar}\bar{\G}$
generated by $\bar{\G}_{\pm}$. These are Hopf subalgebras of
$U_{\hbar}\bar{\G}$, dual to each other if $U_{\hbar}\bar{\G}_{+}$ is
given the coproduct opposite to $\bar\Delta$. 
The pairing $\langle, \rangle_{U_{\hbar}\bar{\G}}$
between $U_{\hbar}\bar{\G}_{-}$ and $U_{\hbar}\bar{\G}_{+}$
is defined by the formulas 
$$
\langle h^{-}[\la], h^{+}[r]\rangle_{U_{\hbar}\bar{\G}} 
={2\over \hbar}\langle r, \la\rangle_{k}, \quad
\langle e[\eps], f[\eta]\rangle_{U_{\hbar}\bar{\G}} 
={1\over \hbar}\langle \eps, \eta\rangle_{k}, 
$$
for $\eps,\eta\in k$, $r\in R,\la\in\La$,
$$
\langle D,K\rangle_{U_{\hbar}\bar{\G}} = 1, \quad 
\langle D, \bar\G(\cK) \rangle_{U_{\hbar}\bar{\G}}=\langle K,
\bar\G(\cK)\rangle_{U_{\hbar}\bar{\G}}=0. 
$$

\begin{prop} \label{??10} (see \cite{ER:qH}, Prop. 6.1) 
The Hopf algebras 
$U_{\hbar}\G$ and $U_{\hbar}\bar{\G}$
 are quasi-triangular, with
respective universal $R$-matrices
\begin{align} \label{form-R}
 & \cR=q^{D\otimes K} q^{ {1 \over 2} \sum_{i\in\NN}
h^{+}[e^{i}]\otimes 
h^{-}[e_{i}] } F, 
\\ & \nonumber 
\bar \cR = F^{21}
q^{ D\otimes K} 
q^{ {1 \over 2} \sum_{i\in\NN}
h^{+}[e^{i}]\otimes h^{-}[e_{i}]} . 
\end{align}
\end{prop}

In fact, $\cR$ and $\bar\cR$ represent the identity for the pairings
$\langle , \rangle_{U_{\hbar}\G}$and $\langle ,
\rangle_{U_{\hbar}\bar\G}$. 

Note that 
$$
\bar \cR= F^{21} \cR F^{-1}  =  e^{\hbar \sum_{i\in\ZZ} f[\eps^{i}] \otimes
e[\eps_{i}] } e^{{\hbar \over 2} D\otimes K}
e^{ {\hbar\over 2} \sum_{i\in\NN}h^{+}[e^{i}]\otimes
h^{-}[e_{i}]} .
$$

\subsection{Regular subalgebra $U_{\hbar}\G_R$} 

Let $U_{\hbar}\G_R$ be the subalgebra of $U_{\hbar}\G$ generated
by the $x[r]$, $x = e,f,h$, $r\in R$.
Then: 
\begin{prop} (see \cite{ER:qH}, sect. 5.2) The inclusion of
$U_{\hbar}\G_R$ in $U_{\hbar}\G$b is a flat deformation of that of 
$U\G_R$ in $U\G$. 
\end{prop}

Moreover, $U_{\hbar}\G_R$ has the following coideal properties
with respect to $\Delta$ and $\bar\Delta$: 
\begin{equation} \label{coideal}
\Delta(U_{\hbar}\G_R) \subset U_{\hbar}\G\otimes 
U_{\hbar}\G_R , 
\quad
\bar\Delta(U_{\hbar}\G_R) \subset U_{\hbar}\G_R 
\otimes U_{\hbar}\G. 
\end{equation}
(see \cite{ER:qH}, Prop. 5.4). 

\subsection{Decomposition of $F$}

Now we would like to decompose $F$ defined in
(\ref{F}) as a product
\begin{equation} \label{decomp}
F_{2}F_{1},\quad \on{ with} \quad  F_{1}\in U_{\hbar}\G\hat\otimes
U_{\hbar}\G_{R}, \quad 
F_{2}\in U_{\hbar}\G_{R}\hat\otimes U_{\hbar}\G.
\end{equation}
The interest of this decomposition lies in the following proposition:  

\begin{prop} \label{??11}
For any decomposition (\ref{decomp}), the map $\Ad
(F_{1})\circ\Delta$ defines an algebra morphism from $U_{\hbar}\G_{R}$ to
$U_{\hbar}\G_{R}\hat\otimes U_{\hbar}\G_{R}$ (where the tensor product
is completed over $\CC[[\hbar]]$). 
\end{prop}
This follows at once from the coideal properties (\ref{coideal}). 

Let us now try to decompose $F$ according to (\ref{decomp}). Let
$(m_{i})$, resp. $(m'_{i})$ be a basis of $U_{\hbar}\G$ as a left,
resp. right $U_{\hbar}\G_{R}$-module. Assume $m_{0}=m'_{0}=1$. 
Due to the form of $F_{1}$ and $F_{2}$, we have decompositions 
$$
F_{2}=\sum_{i}(1 \otimes m'_{j}) F^{(j)}_{2}, \quad
F_{1}=\sum_{i}F^{(i)}_{1}(m_{i} \otimes 1), \quad 
F^{(i)}_{1} , F^{(j)}_{2}\in U_{\hbar}\G_{R}^{\hat\otimes 2}.  
$$
It follows that we have 
\begin{equation} \label{key}
F=\sum_{i}F_{2}F_{1}^{(i)}(m_{i} \otimes 1)
=
\sum_{j}( 1 \otimes m'_{j})F_{2}^{(j)}F_{1}. 
\end{equation}
Let now $\Pi$, resp. $\Pi'$ be the left, resp. right
$U_{\hbar}\G_{R}$-module morphisms from $U_{\hbar}\G$ to
$U_{\hbar}\G_{R}$, such that $\Pi(m_{i})=0$ for $i\neq 0$, $\Pi(1)=1$,
and $\Pi'(m'_{i})=0$ for $i\neq 0$, $\Pi'(1)=1$. 

{}From (\ref{key}) follows that we should have 
\begin{equation} \label{possible}
F_{2}F_{1}^{(0)}=(\Pi\otimes 1)F, \quad 
F_{2}^{(0)}F_{1}=(1 \otimes \Pi')F. 
\end{equation}

Equation (\ref{possible}) determines the possible values of 
$F_{1}$ and $F_{2}$, up to 
right, resp. left multiplication by elements of
$U_{\hbar}\G_{R}^{\hat\otimes 2}$. 

\begin{prop} (see \cite{ER:qH}, Prop. 7.2) 
Let $F_{\Pi,\Pi'} = [(\Pi\otimes 1)F]^{-1} F [(1 \otimes
\Pi')F]^{-1}$; then 
\begin{equation} \label{main}
F_{\Pi,\Pi'} \in U_{\hbar}\G_{R}^{\hat\otimes 2}. 
\end{equation}
\end{prop}

This is the key point of the construction of \cite{ER:qH}; the proof
uses Hopf duality arguments. We prove that 
$F^{-1}[(\Pi\otimes 1)F]$ belongs to 
$U_{\hbar}\G \otimes U_{\hbar}\G_R$, and that 
$[(1 \otimes\Pi')F] F^{-1}$ belongs to 
$U_{\hbar}\G_R \otimes U_{\hbar}\G$. For that, the idea
is to compute the annihilators of the nilpotent parts of 
$U_{\hbar}\G_R$ for the pairings $\langle , \rangle_{U_\hbar\G}$
and $\langle , \rangle_{U_\hbar\bar\G}$; these annihilators are left
and right ideals. Then we pair the second factor of 
$F^{-1}[(\Pi\otimes 1)F]$, and the first factor of 
$[(1 \otimes\Pi')F] F^{-1}$ with this annihilator, and
use the Hopf algebra pairing rules, as well as the algebraic 
properties of $\Pi$ and $\Pi'$, to show that these pairings give
zero. 

{}From the above Prop. follows the solution of the decomposition
problem (\ref{decomp}): 
\begin{prop} \label{??20}
Any decomposition of $F$ according to (\ref{decomp}) is of the form 
$$
F_{2} = [ (\Pi \otimes 1 ) F ] b , \quad 
F_{1} = a [ (1 \otimes \Pi' ) F ], 
$$
with $a,b\in U_{\hbar}\G_{R}^{\hat\otimes 2}$, such that
$ab=F_{\Pi,\Pi'}$. 
\end{prop}

\subsection{Quasi-Hopf structures on $U_{\hbar}\G$ and 
$U_{\hbar}\G_R$} \label{motti}

Let us choose a solution $(F_{1},F_{2})$ of (\ref{decomp}). 
Consider the algebra morphism $\Delta_{R}:
U_{\hbar}\G \to U_{\hbar}\G^{\hat\otimes 2}$, defined as 
\begin{equation} \label{Delta-q-Hopf}
\Delta_{R} = \Ad(F_{1}) \circ \Delta  = \Ad(F_{2}^{-1}) \circ
\bar\Delta;  
\end{equation}
define
\begin{equation} \label{Phi}
\Phi = F_{1}^{23}(1\otimes \Delta)(F_{1}) [F_{1}^{12}(\Delta\otimes
1)(F_{1}) ]^{-1}. 
\end{equation}

\begin{prop} \label{Phi-in-R}
$\Phi$ belongs to $U_{\hbar}\G_{R}^{\hat\otimes 3}$
\end{prop}

Let 
\begin{equation} \label{uR}
u_{R} = m(1\otimes
S)(F_{1}),
\end{equation} 
where $m$ the multiplication of $U_{\hbar}\G$, and 
$S$ is the antipode of $(U_{\hbar}\G,\Delta)$.  

\begin{thm}
The algebra $U_{\hbar}\G$, endowed with the coproduct $\Delta_{R}$,
associator $\Phi$, counit $\varepsilon$, antipode $S_{R} = \Ad(u_{R})
\circ S$, respectively
defined in 
(\ref{Delta-q-Hopf}), (\ref{Phi}), (\ref{counit}), 
(\ref{uR}), and $R$-matrix
$$
\cR_{R} =  [ a^{21}(\Pi'\otimes 1)(F^{21}) ] 
q^{D\otimes K} q^{{1\over 2}\sum_{i\in \NN}h^{+}[e^{i}]\otimes
h^{-}[e_{i}]}
[(\Pi \otimes  1)(F)F_{\Pi,\Pi'}a^{-1}], 
$$
is a quasi-triangular quasi-Hopf algebra. $U_{\hbar}\G_{R}$ is a
sub-quasi-Hopf algebra of it. 
\end{thm}

\subsection{Properties of $q(z,w)$}

In this section, we determine the location of zeroes and poles 
of the function $q(z,w)$, that are at the vicinity of the diagonal. 

First recall that for $f\in \cK$, the product $(f\otimes 1 - 1 \otimes
f)G$ belongs to $\oplus_{s,t\in S}\CC((z_{s}, w_{t}))$. 
Here $z_{s}$ is a local coordinate at $s\in S$, and 
$\CC((z_{s}, w_{t}))
=\CC[[z_{s}, w_{t}]][z_{s}^{-1} , w_{t}^{-1}]$. 

\begin{prop} \label{simply:laced}
Let $z$ be the element of $\cK$ defined as
$(z_{s})_{s\in S}$. Then for some $i \in 1 + \hbar
\oplus_{s,t\in S}\CC((z_{s}, w_{t})) [[\hbar]]$, we have 
$$
z - q^{-\pa} w = i \cdot [
z- w
+ (z- w)G \psi(\hbar , 
(\pa^{i}\otimes 1)\gamma)] . 
$$
\end{prop}

{\em Proof.} Let $\al = (z-w)G$; $\al$ belongs to $\oplus_{s,t\in
S}\CC((z_{s},w_{t}))$. Let us first show that if we replace $z$ by $q^{-\pa}w$,
the expression $ z- w
+ \al(z,w) \psi(\hbar , 
(\pa^{i}_{z}\gamma(z,w))$ vanishes. 

The result of this substitution is a formal series $u(w,\hbar)\in
\cK[[\hbar]]$. It satisfies the equation
\begin{align} \label{toto}
& {{\pa u}\over {\pa\hbar}} (\hbar, w) = 
[-\pa z + (-\pa_{z}\al(z,w))\psi(\hbar,z,w)
\\
&   \nonumber
- \al(z,w)\pa_{z} \psi(\hbar,z,w) + \al(z,w)
{{\pa\psi}\over{\pa\hbar}}(\hbar,z,w)]_{|z =q^{-\pa}w}. 
\end{align}
Since 
$$
\pa_{z}\al(z,w) = (z-w)\gamma(z,w) + [(\pa z) + \al(z,w)] G, 
$$
we have $m(\al) + \pa z =0$ ($m$ being the multiplication map of 
$\cK$); 
it follows that the equation $(z-w)\xi = \pa z + \al(z,w)$ has a solution in
$\oplus_{s,t\in S} \CC((z_{s},w_{t}))$, that we denote by 
$$
{{\pa z + \al(z,w)}\over {z-w}}. 
$$

The l.h.s. of (\ref{toto}) is now equal to 
$$
\left( 
- (\pa z + \alpha(z,w)) (1+ G\psi)(z,w) - 
(z-w + (\al\psi)(\hbar,z,w))(\gamma \psi)(\hbar,z,w) \right)_{|z =
q^{-\pa} w}, 
$$
that is 
$$
-\left((\gamma \psi)(\hbar,z,w) + {{\pa z + \al(z,w)}\over {z-w}}
\right)_{|z = q^{-\pa}w} u(\hbar , w) . 
$$
It follows that the series $u(\hbar , w)$ satisfies the equation
\begin{equation} \label{diff-eq}
{{\pa u}\over {\pa\hbar}} (\hbar, w) =
v(\hbar,w) u(\hbar,w), 
\end{equation}
where $v$ is equal to 
$-\left((\gamma \psi)(\hbar,z,w) + {{\pa z + \al(z,w)}\over {z-w}}
\right)_{|z = q^{-\pa}w} $, and so belongs to $\cK[[\hbar]]$. 
Since we have $u(0,w) =0$, (\ref{diff-eq}) implies that $u$ is
identically zero. 

Let us now recall Lemma 4.2 of \cite{ER:qH}:
\begin{lemma} \label{aux}
Let $z-w+E$ belong to $z-w+\hbar \CC((z,w))[[\hbar]]$, then there exist
unique $e\in \hbar \CC((w))[[\hbar]]$ and $\kappa_{E}\in 1
+\hbar\CC((z,w))[[\hbar]]$, such that 
\begin{equation} \label{arieh}
z-w+E = \kappa_{E} (z-w+e)
\end{equation}
\end{lemma}

Consider the case where $z-w+E = z - w + \al(z,w)\psi(\hbar,z,w)$; the
fact that $u=0$ implies that
$e$ should be such that $z-w+e(w) = z-q^{-\pa}w$; Lemma
\ref{aux} then proves the proposition. 
\qed\medskip

{}From Prop. \ref{simply:laced} follows: 

\begin{corollary} \label{van:q(z,w)}
The function $q(z,w)$ defined by (\ref{q(z,w)}) vansishes for 
$z = q^{-\pa}w$, and its inverse vanishes for $z=q^{\pa}w$. 
\end{corollary}

\section{Quantum currents and quasi-Hopf algebras at level zero} 
\label{level:0}

In this section, we will show how one may define a large family of 
algebras quantizing the ``non-doubly extended'' 
(i.e., original) quasi-Lie bialgebra structures. 
These algebras are defined by relations
similar to those of Thm. \ref{def:quantum:algebra}, using 
functions $q(z,w)$
not necessarily having the zeroes and poles structure described 
by Cor. \ref{van:q(z,w)}. 

Let again $\cK = R \oplus \La$ be some decomposition 
of $\cK$ in isotropic subspaces, and $G = \sum_{i\ge 0} e^i \otimes
e_i$, $e^i,e_i$ dual bases of $R$ and $\La$. We will also set 
$\eps_i = e_i,\eps_{-i-1} = e^i$ for $i\ge 0$. 

Let $a(z,w)$ and $b(z,w)$ belong to 
$1 + \hbar(R\otimes R)[[\hbar]]$ and 
$\hbar(R\otimes R)[[\hbar]]$, and 
define the algebra $U_{a,b}\G$ as the 
algebra generated by the $e[\eps_i],f[\eps_i],h[\eps_i]$, 
with generating series 
$$
e(z) = \sum_{i\in\ZZ} e[\eps_i]\eps_{-i-1}(z), \quad
f(z) = \sum_{i\in\ZZ} f[\eps_i]\eps_{-i-1}(z), 
$$ 
$$
h^+(z) = \sum_{i\ge 0} h[e^i]e_i(z), \quad
h^-(z) = \sum_{i\ge 0} h[e_i]e^i(z), 
$$ 
and the relations
\begin{equation} \label{h:h:level0}
[h[\eps_i] , h[\eps_j]] = 0, 
\end{equation}
for any $i,j$, 
\begin{equation} \label{K+:e:level0}
K^+(z) e(w) K^+(z)^{-1} = {{a+bG^{(21)}} \over {a^{(21)}-b^{(21)}G^{(21)}}}(z,w) 
e(w), 
\end{equation}
\begin{equation} \label{K-:e:level0}
K^-(z) e(w) K^-(z)^{-1} = {{a^{(21)}+b^{(21)}G}\over{a-bG}}(z,w) e(w), 
\end{equation}
\begin{equation} \label{K+:f:level0}
K^+(z)f(w)K^+(z)^{-1}={{a^{(21)}-b^{(21)}G^{(21)}}\over 
{a+bG^{(21)}}}(z,w) f(w), 
\end{equation}
\begin{equation} \label{K-:f:level0}
K^-(z) f(w) K^-(z)^{-1} = {{a-bG} \over 
{a^{(21)}+b^{(21)}G}}(z,w) f(w), 
\end{equation}
\begin{align} \label{e:e:level0}
& (\al(z) - \al(w)) (a(w,z) + b(w,z) G(z,w) ) e(z)e(w)
\\ & = \nonumber
(\al(z) - \al(w)) (a(z,w) - b(z,w) G(z,w) ) e(w)e(z) ,   
\end{align}
\begin{align} \label{f:f:level0}
& (\al(z) - \al(w)) (a(z,w) - b(z,w) G(z,w) ) 
f(z)f(w)
\\ & = \nonumber
(\al(z) - \al(w)) (a(w,z) + b(w,z) G(z,w) )  
f(w)f(z) ,   
\end{align}
for any element $\al$ of $\cK$, 
\begin{equation} \label{e:f:level0}
[e(z) , f(w)] = \delta(z,w) \left( K^+(z) - K^-(z)^{-1} \right) , 
\end{equation}
where $\delta(z,w) = G(z,w) + G(w,z)$, and 
$$
K^+(z) = \exp (\sum_i h[e^i] (1+V)e_i(z)), \quad  
K^-(z) = \exp (\sum_i h[e_i] e^i(z)),
$$ 
and $V$ is the linear operator from $\La$ to $R$ defined 
as follows: let $B$ the linear map from $\La$ to $\cK$ defined
by 
$$
B(\la) = \langle \log {{a+bG^{(21)}}\over{a^{(21)}-b^{(21)}G^{(21)}}}, 
id \otimes \la \rangle ,  
$$ 
and set $B = B_R + B_\La$, where $B_R$ and $B_\La$ are the compositions 
of $B$ with the projections on $R$ and $\La$. Then we have 
$B_\La = \hbar id_{\La} + o(\hbar)$, $B_R = O(\hbar)$, so that
$B_\La$ is invertible, and we set $V = B_R \circ B_\La^{-1}$.  
In other words, if we set 
$$
\log {{a+bG^{(21)}} \over {a^{(21)}-b^{(21)}G^{(21)}}} = \sum_i A_i \otimes e^i , 
$$
we have 
$$
V( (A_j)_{\La}) = (A_j)_R.   
$$
Therefore the relations (\ref{K+:e:level0}), (\ref{K+:f:level0}) 
have correct functional properties with respect to $z$ and 
can be written as
$$
[h[e^i] , e(z)] = B_{\La}(e^i)(z) e(z), \quad
[h[e^i] , f(z)] = - B_{\La}(e^i)(z) f(z), 
$$
for $r\in R$. 

The algebra $U_{a,b}\G$ has coproducts 
$\Delta_{a,b}$ and $\bar\Delta_{a,b}$ defined
by 
\begin{equation} \label{copdt:1}
\Delta_{a,b}(h[\eps_i])  = h[\eps_i] \otimes 1 + 1 \otimes h[\eps_i] ,   
\end{equation}
\begin{align} \label{copdt:2}
\Delta_{a,b}(e(z)) = e(z) \otimes K^+(z) + 1 \otimes e(z),  
\Delta_{a,b}(f(z)) = f(z) \otimes 1 + K^-(z)^{-1} \otimes f(z),  
\end{align}
and 
\begin{equation}
\bar\Delta_{a,b}(h[\eps_i])  = h[\eps_i] \otimes 1 + 1 \otimes h[\eps_i] ,   
\end{equation}
\begin{align}
\bar\Delta_{a,b}(e(z)) = e(z) \otimes 1 + K^-(z)^{-1} \otimes e(z), 
\bar\Delta_{a,b}(f(z)) = f(z) \otimes K^+(z) + 1 \otimes f(z).   
\end{align}

The algebra $U_{a,b}\G$ shares all the properties of 
$U_\hbar\G$ that were used in \cite{ER:qH} for constructing
a quasi-Hopf structure on a subalgebra if it: 

\begin{thm} \label{level:0:qH}
  $(U_{a,b}\G , \Delta_{a,b})$ is a Hopf algebra; it is a flat
  deformation of the enveloping algebra of $\G$. It is the double of
  its subalgebras $U_{a,b}\G_{\pm}$, generated by the modes of
  $K^+(z),e(z)$, respectively $K^-(z),f(z)$.  the Hopf pairing between
  these algebras is defined by
$$
\langle h[e^i],h[e_{-j-1}] \rangle 
= \langle e^i, B_{\La} e_i\rangle_{\cK}, i,j\ge 0,\quad
\langle e[\eps_i], f[\eps_{j}]\rangle = \delta_{i,-j-1} i,j\in \ZZ, 
$$
and the universal $R$-matrix of 
$(U_{a,b}\G , \Delta_{a,b})$ is then equal to
$$
\cR = \exp\left({1\over 2}\sum_{i\ge 0} h[e^i] \otimes
h[B_{\La} e_i]\right) F_{a,b},  
$$ 
where $F_{a,b} = \sum_i \al^i \otimes \al_i$, $(\al^i),(\al_i)$ dual
bases of the subalgebras $U_{a,b}\N_+$ and $U_{a,b}\N_-$ of
$U_{a,b}\G_+$ and $U_{a,b}\G_-$, respectively generated by the
$e[\eps],\eps\in\cK$ and the $f[\eps],\eps\in\cK$. 

The same statements hold with $\Delta_{a,b}$ replaced by 
$\bar\Delta_{a,b}$. 
The subalgebras are then generated by the modes of $K^-(z),e(z)$, 
respectively $K^+(z),f(z)$. The Hopf pairing is defined by the 
same formulas as above, and the universal $R$-matrix is equal to 
$$
\bar\cR = F^{21}_{a,b}
\exp\left({1\over 2}\sum_{i\ge 0} h[e^i] \otimes h[B_{\La} 
e_i]\right) . 
$$ 
Moreover, $\bar\Delta_{a,b}$ is obtained from $\Delta_{a,b}$ by the 
twist $\bar\Delta_{a,b} = F_{a,b}\Delta_{a,b} F_{a,b}^{-1}$. 

The subalgebra $U_{a,b}\G_R$ spanned by the $x[e^i]$, $x=e,f,h$, $i\ge 0$ 
is a flat deformation of the enveloping algebra of $\G_R = \bar\G \otimes R$. 
It satisfies 
$$
\Delta_{a,b}(U_{a,b}\G_R) \subset U_{a,b}\G \otimes U_{a,b}\G_R, 
\quad 
\bar\Delta_{a,b}(U_{a,b}\G_R) \subset U_{a,b}\G_R \otimes U_{a,b}\G. 
$$
\end{thm}

The procedure of \cite{ER:qH} can then be followed to
obtain a decomposition of $F$, which will serve to define
a quasi-Hopf algebra structure on $U_{a,b}\G_R$.   

\begin{remark} {\it Dependence of $U_{a,b}\G$ in $(a,b)$.}
\label{isoms}
Clearly, $U_{a,b}\G$ depends on $(a,b)$ only through 
$$
q_{a,b}(z,w) = \left({{a+bG^{(21)}} \over {a^{(21)}-b^{(21)}G^{(21)}}} 
\right)(z,w); 
$$
in particular, we have $q_{a,b} = q_{sa,sb}$, if  
$s\in (R\otimes R)[[\hbar]]$ is invertible and 
symmetric in $(z,w)$, and also  
$q_{a,b} = q_{a - c G,b + c}$, if $c$ is a function of $R\otimes R$ 
vanishing on the diagonal of $X$.

On the other hand, if we have $q_{a,b} = \la q_{a',b'}$ for some 
invertible $\la$ in $(R\otimes R)[[\hbar]]$, there is an isomorphism 
between the algebras $U_{a,b}\G$ and $U_{a',b'}\G$, defined by 
multiplying the fields $e(z)$ and $f(z)$ by suitable combinations of 
$K^+(z)$ and $K^-(z)$. 
\hfill \qed \medskip
\end{remark}

\begin{remark} The condition that $U_{a,b}\G$ can be 
embedded in some algebra ``with derivation and central 
extension'' quantizing its double extension seems 
to impose severe constraints to $a$ and $b$. Indeed, a 
natural way to achieve this is to use the relations of 
\cite{ER:qH}, sect. 8. These relations imply in particular 
that we have 
$$
K^+(z)K^-(w)K^+(z)^{-1}K^-(w)^{-1} = 
{{q_{a,b}(z,w)} \over {q_{a,b}(z,q^{K\pa}w)}} , 
$$ 
$$
[D, K^+(z)] = - \pa_z K^+(z) + (Ah^+)(z) K^+(z),
$$
$$
[D, K^-(z)] = - \pa_z K^-(z) + (Bh^+)(z) K^-(z),  
$$
where $A$ and $B$ are some finite rank operators from 
$\La$ to $R$, and $K$ is the central generator; 
this implies in particular that the set of zeroes and poles
of $q(z,w)$ is stable by the diagonal action of $\pa$ on 
$X \times X$. In the case studied in \cite{ER:qH}, 
these sets are $\{ (x,q^{\pa}x), x\in X\}$ and the diagonal of 
$X \times X$. 
\end{remark}

\section{Examples on a rational curve} \label{rat:sect}

\subsection{Manin pairs}
In this section, we will consider the following situation.
Let us fix an integer $N\ge 2$. Let us set $\cK = \CC((z))$, 
$\omega = z^{N-1}dz$.  

\subsubsection{} \label{even}
If $N$ is odd, write $N = 2n+1$, with $n$ an integer
$\ge 0$. Let us set 
$R = z^{-n-1}\CC[z^{-1}]$, $\La = z^{-n}\CC[[z]]$.
Then $R$ is a maximal isotropic subring of $\cK$ for the pairing
induced by $\omega$, and $\La$ is a maximal isotropic supplementary. 

Dual bases of $R$ and $\La$ are 
$e^{i} = z^{-n-i-1}$, $e_{i} = z^{i-n}$ for $i \ge 0$. 
We then have 
$$
G = \sum_{i\ge 0} z^{-n-i-1} \otimes z^{i-n}
= 
{{(zw)^{-n}} \over{z-w}} , 
$$
expanded for $w$ near $0$. 

We construct then a Manin pair as follows: we define $\G$ as the 
Lie algebra $(\bar\G \otimes\cK) \oplus \CC K \oplus \CC D$, 
where the central and cocentral extensions are defined as in sect. 
\ref{def:sect}, endowed the usual scalar product. The Lie subalgebra  
$\G_R = (\bar\G\otimes R) \oplus \CC K$ is then maximal 
isotropic for this scalar product. This defines a Manin pair. 
Quasi-Lie bialgebra structures are then defined on $\G$ and on 
$\G_R$ by the choice of the isotropic complement 
$\G_{\La} = (\bar\G\otimes \La) \oplus \CC D$ of $\G_R$. 

\subsubsection{} \label{odd}

Suppose $N$ is even. Write $N = 2(n+1)$, with $n$ an integer 
$\ge 0$. 

Let $\wt\G$ be the semidirect product $\wt\G = \G \oplus
\CC \wt h$, where $\G\subset\wt\G$ is a Lie algebras embedding, and 
the action of $\wt h$ on $\G$ is such that $\wt h - h\otimes z^{-n}$
is central. Extend the scalar product $\langle , \rangle_{\G}$ 
of $\G$ to a scalar product $\langle , \rangle_{\wt\G}$  on $\wt\G$ 
by the rules $\langle \wt h,\G \rangle_{\wt\G} = 0$, 
$\langle \wt h,\wt h \rangle_{\wt\G}  + \langle h\otimes z^{-n},
h\otimes z^{-n} \rangle_{\wt\G}=0$.  

A quasi-Lie bialgebra structure on $\wt\G$ is then defined as follows: 
let $\G_R^{(0)}$ and $\G_\La^{(0)}$ be the subspaces of $\G$ equal to 
$(\bar\G\otimes z^{-n-1}\CC[z^{-1}])\oplus\CC K$ and 
$(\bar\G\otimes z^{1-n}\CC[[z]])\oplus\CC D$, 
and define $\wt\G_R$ and $\wt\G_\La$ as the direct sums of 
$\CC(\wt h - h\otimes z^{-n})$ and $\CC(\wt h + h\otimes z^{-n})$ 
with the images of $\G_R^{(0)}$ and $\G_\La^{(0)}$ in $\wt\G$. 
Then $\wt\G_R$ is a maximal
isotropic Lie subalgebra of $\wt\G$ and $\wt\G_\La$ is an isotropic 
complement. This defines Lie quasi-bialgebra structures on $\wt\G$ 
and $\wt\G_R$. 

By the natural projection of $\wt\G$ to $\G$, these structures
define quasi-Lie bialgebra structures on $\G$ and on $\G_R
= (\bar\G\otimes z^{-n}\CC[z^{-1}])\oplus \CC K$; the
structure on $\G$ is not a double one. 

\begin{notation}
Here and later, we will use the notation 
$z_{\la} = (z^N + \la N \hbar)^{1/N}$. We have 
$z_{\la} = q^{\la\pa}(z)$, where $\pa$ is the derivation defined 
by $\omega$. 
\end{notation}

\subsection{The functions $q(z,w)$} 

Define the series $\phi(z,w)$ of $\CC[[z^{-1}]]((w))[[\hbar]]$
as the expansion of $\log{{z_1-w}\over{z-w}}$. Set $\phi(z,w)
= \sum_{p,q\in\ZZ} a_{pq}z^{p}w^{q}$, and 
$\phi_{w^{>-n}}(z,w)=\sum_{q>-n,p\in\ZZ}a_{pq}z^{p}w^{q}$. 
It is easy to check that $\phi_{w^{> -n}}(z,w)$ belongs to 
$$
z^{-n}w^{-1}\CC[[z^{-1},w^{-1},\hbar]].
$$

\begin{prop} \label{nesher}
1) Let the notation be as in \ref{even}.  
Set $N = 2n+1$. There exists some linear 
operator $U:\La\to \hbar R[[\hbar]]$ such that 
\begin{equation} \label{hanukah}
(1\otimes ( {{q^{\pa} - q^{-\pa}}\over{\pa}} + U)) G 
=
\log {{ z_1 - w}\over{{z - w_1}}} 
+\phi_{z^{> -n}}(z,w) - \phi_{w^{>-n}}(w,z)
. 
\end{equation}

2) Let the notation be as in \ref{odd}. Let $N = 2n+2$. 
For some linear operator $U: z^{-n}\CC[[z]] \to 
z^{-n}\CC[z^{-1}][[\hbar]]$, (\ref{hanukah}) holds. 

\end{prop}

{\em Proof.} Let us prove 1). 
Let us first show that 
\begin{equation} \label{differ}
D_{\hbar} = G - (q^{-\pa}\otimes q^{-\pa}) G 
\end{equation}
belongs to $(R\otimes R)[[\hbar]]$. 

We have 
\begin{equation} \label{aigle}
(\pa\otimes 1 + 1\otimes \pa)G \in R\otimes R
\end{equation}
This follows from Lemma \ref{barbara}. 
To show (\ref{aigle}), we may also compute explicitly 
$$
(\pa\otimes 1)G = -G^2 + \gamma, 
$$
with $\gamma = -z^{-2n}w^{-n}{1\over{z-w}}
[{{z^{-n}-w^{-n}}\over{z-w}} + nz^{-n-1}]$, which belongs to 
$R\otimes R$ because the term in brackets vanishes for $z=w$, 
and 
$$
(1\otimes \pa)G = G^2 - \gamma^{(21)}, 
$$
so that 
$(\pa\otimes 1 + 1\otimes \pa)G = \gamma - \gamma^{(21)}$
belongs to $R\otimes R$. 

Now $R$ is stable under $\pa$, so that 
$D_{\hbar} = [(q^{-\pa}\otimes q^{-\pa}) - 1]G = 
{{q^{-(\pa\otimes 1+1\otimes\pa)}-1}\over
{\pa\otimes 1+1\otimes\pa}}
(\pa\otimes 1 + 1\otimes \pa)G$ also belongs to 
$R\otimes R$. 

Therefore $(q^{\pa}\otimes 1)(D_\hbar)$
also belongs to $(R\otimes R)[[\hbar]]$. It follows that 
for some linear operator $V_+: \La \to R[[\hbar]]$, we have  
$$
(q^{\pa}\otimes 1)(D_\hbar) = (1\otimes V_{+})(G). 
$$
Therefore
$$
(1\otimes (q^{-\pa} + V_{+})) G =
(q^{\pa}\otimes 1) \left( 
{{z^{-n}w^{-n}}\over{z-w}} \right) 
= 
{ {z_1^{1-N} }
\over { z_1 - w }} 
- \left( 
{ { z_1^{1-N} } 
\over { z_1 - w }} 
\right)_{w^{>-n}} 
; 
$$
let $U_+$ be the unique linear operator from $\La$ to 
$\hbar R[[\hbar]]$,
such that $\pa_{\hbar} U_+ = V_{+}$. Integrating in $\hbar$, 
we obtain 
\begin{equation}
(1\otimes {{1-q^{-\pa}}\over{\pa}} + U_{+})) G 
=
\log {{z_1 - w}\over{{z - w}}}
-\phi_{w^{>-n}}(z,w)
. 
\end{equation}

We may construct in the same way a linear operator $U_-$
from $\La$ to $\hbar R[[\hbar]]$, such that 
\begin{equation}
(1\otimes ( {{q^{\pa}-1}\over{\pa}} + U_{-})) G 
=
\log { {z-w}\over {z-w_1}} +\phi_{z^{> -n}}(w,z). 
\end{equation}
To obtain the statement of the proposition, we then set 
$U = U_+ + U_-$. 

The proof of 2) is similar. 
\hfill \qed \medskip 

Let us choose $U$ like in Prop. \ref{nesher}. We then 
find 
\begin{equation} \label{q:z^n}
q(z,w) = \exp( \phi_{z^{>-n}}(z,w) - 
\phi_{w^{>-n}}(w,z) ) 
{{ (z^N + \hbar N)^{1/N} - w}\over
{z - (w^N + \hbar N)^{1/N}}}. 
\end{equation}

\subsection{The algebra $U_{\hbar,z^{N-1} dz}\G$}  

We denote by $U_{\hbar,z^{N-1}dz}\G$ the Hopf algebra 
resulting from the construction of Thm. \ref{def:quantum:algebra}. 
It contains a regular subalgebra, generated by the 
$x_i$, $i\leq -n$. 

In what follows, we will set
$$
U_{\hbar,z^{-n-2}dz}\G = U_{\hbar,z^{n}dz}\G 
$$ 
for $n\ge 0$, and $U_{\hbar,z^{-1}dz}\G$ equal to the quantum 
affine algebra
attached to $\G$. 

One interest of the algebras $U_{\hbar,z^{N-1} dz}\G$ lies in 
the following

\begin{thm} \label{isomorphisms}
Let $U_{\hbar}\G$ be the algebra of Thm. \ref{def:quantum:algebra}, 
attached to the data $(X,\omega,S)$ and let $U_{\hbar}\G'$ be 
its subalgebra with the same generators except $D$. Let for 
each point $s$ of $S$, $n_s$ be the order of the zero or pole 
of $\omega$ at $s$. Then $U_{\hbar}\G'$ is isomorphic to the 
quotient $\otimes_{s\in S} U_{\hbar,z^{n_s}dz}\G' / (K^{(s)} 
- K^{(t)})$, 
where $K^{(s)}$ is the central generator of the $s$th factor. 
\end{thm}

{\em Proof.} The argument is similar to that of 
\cite{coinvariants}, introd.: 
fix at each point $s$ a coordinate $z_s$ such that $\omega$ is 
locally expressed by $z_s^{n_s}dz_s$. Then for each $s$, we have 
a specialization morphism $\ev_s$ from $U_{\hbar}\G'$ to 
$U_{\hbar,z^{n_s} dz}\G'$, sending each $x[\eps_t]$ to 
$\delta_{st}x[\eps_t]$ if $\eps_t\in \cK_t$ and $K$ to $K$. Fix a 
coproduct $\Delta_R$ for $U_{\hbar}\G'$ as in sect. \ref{motti}. 
We have then an algebra morphism $\Delta_R^{(\card S)}$ from 
$U_{\hbar}\G'$ to $(U_{\hbar}\G')^{\otimes \card S}$, defined by 
$\Delta_R^{(\card S)} = (\Delta_R \otimes id^{\otimes \card S - 1}) 
\circ \cdots \circ \Delta_R$. Choose an order of the points of 
$S$ and compose $\Delta_R^{(\card S)}$ with $\otimes_{s\in S}\ev_s$. 
The resulting map is an algebra morphism from  $U_{\hbar}\G'$ to 
$\otimes_{s\in S} U_{\hbar,z^{n_s}dz}\G'$. That it gives an 
isomorphism after composition with projection to the quotient 
by the ideal generated by the $K^{(s)} - K^{(t)}$ follows from
inspection of its classical limit. 
\hfill \qed \medskip 

The algebras $U_{\hbar,z^{N-1}dz}\G$ have also the 
property that for any nonzero complex $\la$, 
$U_{\la\hbar,z^{N-1}dz}\G$  is isomorphic with 
$U_{\hbar,z^{N-1}dz}\G$; this follows from the fact that (writing
the formal parameter in indices)  
$q_{\al^N \hbar}(\al z,\al w) = q_{\hbar}(z,w)$. This generalizes
the properties of Yangians of being isomorphic for all nonzero 
values of the deformation parameter. 

\begin{remark}
In the framework of the preceding sections, one 
should consider the curve $X = \CC P^{1}$ with differential 
$\omega$ and marked points $0$ and $\infty$. The resulting algebra 
would be nothing but the tensor square of $U_{\hbar,z^{N-1}dz}\G$. 
\hfill \qed \medskip 
\end{remark}

\begin{remark} \label{completion} 
If we complete $U_{\hbar,z^{N-1}dz}\G$ with respect to the ideals 
generated by the $x[z^{-i}]$, $i\ge n$, $x = e,f,h$, the relations 
defining it make sense for complex values of $\hbar$. 
\end{remark}

\subsection{Another presentation of the vertex relations of 
$U_{\hbar,z^{N-1}dz}\G$}

It is easy to see that after we multiply the generating series 
$e(z)$ by a suitable Cartan fields, they satisfy
$$
((z^N+N\hbar)^{1/N}-w)\wt e(z)\wt e(w) = (z-(w^N+N\hbar)^{1/N}) 
\wt e(w)\wt e(z).  
$$
We will show how this relation can be written 
avoiding the use of $N$th roots.  

Let us denote by $\ZZ_N$ the group $\ZZ/N\ZZ$ and by 
$\mu_N$ the group of $N$th roots of unity in 
$\CC$. Let us decompose the field
$\wt e(z)$ as
\begin{equation} \label{x:y:z^n}
\wt e(z) = \sum_{\al\in\ZZ_{N}} e^{(\al)}(z), \quad \on{with} 
\quad 
e^{(\al)}(\zeta z) = \zeta^\al e^{(\al)}(z),
\end{equation}
for $\zeta\in\mu_N$. 
We also set $e^{(\al)}(z) = z^\al E^{(\al)}(z^n)$, 
where we denote by 
$\bar\al$ 
the representative in $[0,N-1]$ of the element $\al$ of $\ZZ_N$, and
we abuse notations by writing $z^\al = z^{\bar \al}$. 

Define for $a,b\in \ZZ_N$, $r(a,b)$ as the number (equal to $0$ or 
$1$) such that $\bar a + \bar b = \overline{a+b} + r(a,b)N$; this 
is the carry over for the addition of $\bar a$ and $\bar b$ mod. 
$N$.  

\begin{prop} \label{plantes}
The relation (\ref{x:y:z^n})
is equivalent to the system 
of relations 
\begin{align} \label{X:Y}
& (Z - W + N\hbar)
\sum_{\al\in\ZZ_N} 
Z^{r(N-p-\al,\al)}(W+ N\hbar)^{r(p+\al-1,q-\al)} 
E^{(\al)}(Z)
E^{(q-\al)}(W)
\\ & \nonumber =
(Z - W - N\hbar) 
\sum_{\al\in\ZZ_N} 
(Z+N\hbar)^{r(N-p-\al,\al)}W^{r(p+\al-1,q-\al)} 
E^{(q-\al)}(W) E^{(\al)}(Z). 
\end{align}
\end{prop}

{\em Proof.} 
Write the relation (\ref{x:y:z^n}) as 
$$
{{\wt e(z)\wt e(w)}\over{z-w_1}} = {{\wt e(w)\wt e(z)}\over{z_1-w}}.  
$$
It implies that for $p\in\ZZ_N$, 
$$
\sum_{\zeta\in\mu_N} {{\zeta^p \wt e(\zeta z)\wt e(w)}
\over{\zeta z-w_1}} = 
\sum_{\zeta\in\mu_N} {{\zeta^p \wt e(w)\wt e(\zeta z)}
\over{\zeta z_1-w}} . 
$$
Since we have
$$
\sum_{\zeta\in\mu_N} {{\zeta^p}\over{\zeta z -w}} = 
{{n w^{p-1}z^{N-p}}\over{z^N - w^N}}, 
$$
it follows that 
$$
\sum_{\zeta\in\mu_N} {{\zeta^p \wt e(\zeta z)}\over{\zeta z-w}} =
{{N}\over{z^N - w^N}}\left( \sum_{\al\in\ZZ_N} 
z^{N-p-\al}w^{p+\al-1} e^{(\al)}(z)
\right) . 
$$
Therefore 
\begin{align*}
{{N}\over{z^N - w_1^N}}
& \left( \sum_{\al\in\ZZ_N} 
z^{N-p-\al}w_1^{p+\al-1} e^{(\al)}(z)
\right) 
\sum_{\beta\in\ZZ_N} e^{(\beta)}(w)
\\ & =
{{N}\over{z_1^N - w^N}}
\left( \sum_{\beta\in\ZZ_N} e^{(\beta)}(w) \right) 
\left( \sum_{\al\in\ZZ_N} 
z_1^{N-p-\al}w^{p+\al-1} e^{(\al)}(z)
\right).  
\end{align*}
Separating isotypic components for the action of $\ZZ_N$ in 
the variable $w$, we get for each $q\in\ZZ_N$  
\begin{align*}
{{N}\over{z^N - w^N -  N\hbar}}
& \sum_{\al\in\ZZ_N} 
z^{N-p-\al}w_1^{p+\al-1} e^{(\al)}(z)
e^{(q-\al)}(w)
\\ & =
{{N}\over{z^N - w^N +N\hbar}}
\sum_{\al\in\ZZ_N} 
z_1^{N-p-\al}w^{p+\al-1} e^{(q-\al)}(w) e^{(\al)}(z),   
\end{align*}
so that in terms of fields $X^{(\al)}$ we obtain  
\begin{align*}
& {{N}\over{z^N - w^N - N\hbar}}
\sum_{\al\in\ZZ_N} 
z^{N-p-\al}z^{\al}w_1^{p+\al-1}w_1^{q-\al} E^{(\al)}(z^N)
E^{(q-\al)}(w^N)
\\ & =
{{N}\over{z^N - w^N + N\hbar}}
\sum_{\al\in\ZZ_N} 
z_1^{N-p-\al}z_1^{\al}w^{p+\al-1}w^{q-\al} E^{(q-\al)}(w^N) 
E^{(\al)}(z^N),   
\end{align*}
so that we obtain, with 
$Z = z^N,W = w^N$,   
\begin{align*}
& {{N}\over{Z - W -  N\hbar}}
\sum_{\al\in\ZZ_N} 
Z^{r(N-p-\al,\al)}(W+ N\hbar)^{r(p+\al-1,q-\al)} 
E^{(\al)}(Z)
E^{(q-\al)}(W)
\\ & =
{{N}\over{Z - W +  N\hbar}}
\sum_{\al\in\ZZ_N} 
(Z+ N\hbar)^{r(n-p-\al,\al)}W^{r(p+\al-1,q-\al)} 
E^{(q-\al)}(W) E^{(\al)}(Z),  
\end{align*}
that is (\ref{X:Y}).  

The above arguments can easily be reversed to show the 
proposition. \hfill \qed 

\begin{remark}
We may construct an algebra $A_{\hbar,z^{N-1}dz}$ with generators 
$E^{(\al)}_i$, $i\in\ZZ$, arranged in series 
$$
E^{(\al)}(Z) = \sum_{i\in\ZZ} E^{(\al)}_i Z^{-i}, 
$$ 
subject to the above relations (\ref{X:Y}). As we have seen, for
$\hbar$ a formal parameter, it is isomorphic with 
the part of $U_{\hbar,z^{N-1}dz}$ generated by the field 
$\wt e(z)$. This is also true in the case when $\hbar$ is
complex, after
we complete $U_{\hbar,z^{N-1}dz}$ as in Rem. \ref{completion}. 
However, since $\hbar$ appears polynomially 
in the defining relations of $A_{\hbar,z^{N-1}dz}$, they make 
sense without completing the algebra, even when $\hbar$ is complex. 
\end{remark}
\begin{remark} The algebra
$A_{\hbar,z^{N-1}dz}$ has an obvious
morphism to the upper nilpotent subalgebra of the double Yangian
$DY(\SL_2)$, defined by $E^{(\al)}(Z)\mapsto \delta_{\al 0}e(Z)$.  
\end{remark}

\begin{remark} Prop. \ref{plantes} can easily be extended to the 
case of mixed vertex relations 
$$
(z_\la - w) x(z) y(w) = (z-w_\la) y(w) x(z). 
$$
\end{remark}

\section{Genus $>1$ examples associated to odd 
theta-characteristics} \label{g>1}

Let $X$ be a smooth curve of genus $>1$. Let $\omega$ be a regular
form on $X$ all whose zeroes are double. The existence of such a 
form follows from that of a nonsingular odd theta-characteristic -- that 
is, from the existence of an effective divisor with double equivalent to 
the canonical divisor (see e.g. \cite{Mum}, Lemma 1, p. 3.208 -- 
or \cite{Fay}). Let $\delta = \sum_{i=1}^{g-1} \delta_i$ be this effective
divisor; we then have $\ddiv(\omega) = 2 \delta$. Let $\cL_{\delta}$
be the line bundle associated with $\delta$. Then we have 
$\cL_{\delta}^{\otimes 2} = K$ and $h^0(\cL_{\delta}) = 1$. 

Let us also recall the properties of the vector of Riemann constants
(\cite{Fay,Mum}). 
Let $\Jac^n(X)$ be the degree $n$ component of the Jacobian 
of $X$. View the basic theta-function $\theta$ as a quasi-periodic
function on a cover of $\Jac^0(X)$. We denote the same way points of $X$
and their image in $\Jac^1(X)$. Then for some vector $\Delta$ of 
$\Jac^{g-1}(X)$, we have
\begin{equation} \label{Riem}
\theta( - \Delta + \sum_{i=1}^{g-1} y_i) = 0, 
\end{equation}
for any collection of $g-1$ points $y_i$ of $X$. 
Moreover, the zero set of $\theta$ is equal to $\{ - \Delta + 
\sum_{i=1}^{g-1} y_i, y_i\in X\}$.

\subsection{Quasi-Lie bialgebras}

Here we will consider some Manin pairs, where the Lie subalgebra
will be formed by the currents regular at some points, 
as it was done in \cite{ell:QG} in genus $1$. 

More precisely, let $S'$ be a set of points of $X$ not containing any
$\delta_i$, and let us define 
$$
\cK = \oplus_{s\in S'} \cK_s, \quad \cO = \oplus_{s\in S'} \cO_{s}, 
$$
and the pairing $\langle , \rangle_{\cK}$ on $\cK$ by 
$$
\langle f,g \rangle_{\cK} = \sum_{s\in S'} \res_s (fg\omega). 
$$
$\cO$ is clearly an isotropic subalgebra of $\cK$. 

\subsubsection{Isotropic subspaces} \label{suppl:sect}

We define some isotropic subspace $L$ of $\cK$ as follows. 
Fix a system $(a_i,b_i)_{i=1,\cdots,g}$ 
of $a$- and $b$- cycles on $X$. Let us denote by 
$\wt X$ the universal cover of $X$, and by $\gamma_{a_i}$ and 
$\gamma_{b_i}$ the deck transformations associated to 
$a_i$ and $b_i$. Let $X^{(a)}$ 
be the quotient of $\wt X$ by the equivalences $z\sim \gamma_{a_i}z$.
Let us fix lifts $a_i$ of the $a$-cycles in $X^{(a)}$, that we 
also denote by $a_i$, and let $X_0$ be 
the fundamental domain in $X^{(a)}$, bounded by the $a_i$ and the 
$\gamma_{b_i}(a_i)$. We identify local fields at points of 
$S'$ with the local fields at their lifts in $X_0$, and denote by $\wt
S'$ the lift of $S'$ to the fundamental domain.

Define $L$ as the set of expansions at the points of $\wt S'$ of the
functions $f$ such that 
$$
f(\gamma_{a_i}z) = f(z), \quad f(\gamma_{b_i}z) = f(z) + c_i(f), 
$$
$f$ is regular except at the points of $\wt S'$, has simple 
poles at most at the lifts of the $\delta_i$, and is such that
$$
\int_{a_i} f\omega = - c_i(f)/2 \int_{a_i}\omega. 
$$

We can generalize this construction of isotropic subspaces 
of $\cK$ as follows. 
Let $V$ be a vector subspace of 
$\CC^g$. Define $L_V$ as the set of expansions at the points of 
$\wt S'$ of the functions $f$ defined on $X^{(a)}$,
such that $f(\gamma_{a_i}z) = f(z)$, 
$(f(\gamma_{b_i}z) - f(z))_{i=1,\cdots,g}$ belongs to $V$, and 
the periods condition 
\begin{equation} \label{cond:isot:LV}
\sum_i \al_{ij} \left( \int_{a_i} f \omega 
+\int_{\gamma_{b_i}(a_i)} f \omega \right) 
= 0, 
\quad j = 1, \cdots, s,  
\end{equation}
where $(\alpha_{ij})_{i=1,\cdots,g}$, $j = 1,\cdots, s$ 
are the coordinates of a basis of $V$. 

We then have: 

\begin{lemma} For each subspace $V$ of $\CC^g$, 
$L_V$ is isotropic for 
$\langle , \rangle_{\cK}$ 
\end{lemma}

{\em Proof.}
Let $f,g$ belong to $L_V$. By the residues theorem, 
$\langle f,g \rangle_{\cK}$ is equal to 
\begin{equation} \label{solts}
- \sum_{i=1}^{g-1} \res_{\delta_i}(fg\omega)
- \sum_{i=1}^{g} \left( 
\int_{a_i} fg\omega - \int_{\gamma_{b_i}(a_i)} fg\omega
\right) ; 
\end{equation}
by the simple poles conditions on $f$ and $g$, the first
sum in (\ref{solts}) vanishes. On the other hand, set 
\begin{equation} \label{shifts}
f(\gamma_{b_i}z) - f(z) = 
\sum_{j=1}^{s} \la_{ij}(f)\al_{ij}, 
\quad
g(\gamma_{b_i}z) - g(z) = 
\sum_{j=1}^{s} \la_{ij}(g)\al_{ij};  
\end{equation}
we then have 
\begin{equation} \label{james}
\sum_{i=1}^g \al_{ij} \left( 2 \int_{a_i} f \omega
+ \sum_{j=1}^{s} \la_{ij}(f)\al_{ij} 
\int_{\gamma_{b_i}(a_i)}\omega
\right) = 0, 
\end{equation}
\begin{equation} \label{bond}
\sum_{i=1}^g \al_{ij}\left( 2 \int_{a_i} g \omega
+ \sum_{j=1}^{s} \la_{ij}(g)\al_{ij} 
\int_{\gamma_{b_i}(a_i)}\omega
\right) = 0 , 
\end{equation}
$j = 1,\cdots,s$. 
It follows from (\ref{shifts}) that 
the second sum of (\ref{solts}) is equal to 
\begin{equation} \label{tps} 
\sum_{i=1}^{g}  
\int_{a_i} (f + \sum_{j=1}^{s} \la_{ij}(f)\al_{ij}
 ) (g + \sum_{j=1}^{s} \la_{ij}(g)\al_{ij} )\omega 
- \int_{a_i} fg\omega ; 
\end{equation}
multiplying (\ref{james}) by $\la_{ij}(g)$ 
and (\ref{bond}) by $\la_{ij}(f)$ and summing up both 
sets of equations, we find that (\ref{tps}) vanishes.
Therefore (\ref{solts}) vanishes.  
\hfill \qed \medskip

The spaces $L_V$ differ from $L$ only by finite-dimensional pieces 
(that is, their projection to $L$ parallel to $\cO$ has 
finite kernel and cokernel). For $V = \CC^g$, we have 
$L_V = L$.

\begin{remark} {\it Lagrangian supplementaries associated with bundles.} 
Fix a family $(g_i)_{i = 1,\cdots,g}$ of elements of $G$. 
We may consider the subspace $L_{(g_i)}$ of $\G \otimes\cK$ formed
by the expansions at the points of $S'$ of the functions 
$f$ from $\wt X$ to $\G$, such that $f$ has simple poles at $\delta$, 
$$
f(\gamma_{a_i}z) = f(z) , \quad f(\gamma_{b_i}z) = 
\Ad(g_i)(f(z)) + x_i(f),  
$$
and $\int_{a_i} f\omega = - {1\over 2} \Ad(g_i^{-1})x_i(f) 
\int_{a_i}\omega$. The sum of
the residues of $f$ at the points of $\delta$ is then
$\sum_i (1 - \Ad(g_i^{-1})) 
(x_i(f))$ (which needs not be zero).  
$L_{(g_i)}$ is an isotropic subspace of $\G\otimes\cK$. 

In the case of sums of line bundles, we obtain the analogues 
of the spaces $L_{\la}$ of \cite{ell:QG}. 

It would be interesting to understand if the $r$-matrix associated
with these supplementaries satisfies some variant of the 
dynamical Yang-Baxter equation. 
\end{remark}

\subsubsection{Green kernels} \label{green:sect}

In this section, we will consider the case when $V$
is one-dimensional; set $V = \CC h$, $h\in \CC^g$. 

Let us set for $z,w$ in $X^{(a)}$, 
$$
G_h(z,w) = {{\pa_h \theta(z-w+\delta-\Delta)}\over 
{\theta(z-w+\delta-\Delta)}}. 
$$
This function has the following properties: 

\begin{prop}
$G_h(z,w)$ is antisymmetric in $z$ and $w$. It has poles
for the projections on $X$ of $z$ or $w$ equal, or equal 
to some $\delta_i$. We have for any $z,w$,
\begin{equation} \label{q-per}
G_h(\gamma_{b_i} z,w) =  G_h(z,w) - h_i.  
\end{equation}
Near the diagonal $z=w$, $G_h(z,w)$ has the expansions
\begin{equation} \label{expansion}
G_h(z,w) = {C(h)\over{\int_z^w \omega}}+ O(1), 
\end{equation}
for some constant $C(h)$, which is non-zero iff $h$ does not
belong to the linear span of the $V_{\delta_i}$; for any 
point $x$ of $X$, we denote by $V_x$ some tangent vector at $x$ 
of the embedding of $X$ in its Jacobian. 
\end{prop}

{\em Proof.}
The function defined on $\CC^g$ by 
$\zz \mapsto \theta(-\Delta + \delta + \zz)$ is odd, so that
$\zz \mapsto \pa_h\theta/\theta(-\Delta + \delta + \zz)$ is 
also odd; it follows that $G_h(z,w)$ is antisymmetric. 

Recall that 
$$
\theta(\zz + A_i) = const \cdot \theta(\zz), 
\theta(\zz + B_i) = const \cdot e^{-z_i} \theta(\zz), 
$$
where $\zz = (z_i)_{1\leq i\leq g}$, 
$A_i$ are the basis vectors of $\CC^g$ and $B_i$
the vector $(\int_{b_i} \omega_1, \cdots ,$ $\int_{b_i}\omega_g)$, 
and the $\omega_i$ are the holomorphic differentials such that
$\int_{a_i}\omega_j = \delta_{ij}$. 
Taking logarithmic derivative, we find that 
$$
(\pa_h \theta / \theta)(\zz + A_i) = (\pa_h \theta / \theta)(\zz), 
(\pa_h \theta / \theta)(\zz + B_i) = (\pa_h \theta / \theta)(\zz)
- h_i .  
$$
if $h$ has components $(h_1,\cdots,h_g)$. (\ref{q-per}) follows 
from these identities. 

Finally, to prove (\ref{expansion}), we need the following 
result: 

\begin{lemma} \label{diagonal}
The expression 
$\alpha(z) dz = d_z \theta(z-w+\delta - \Delta)_{| z=w}$
is a $1$-form on $X^{(a)}$, defined as the restriction 
on the diagonal of $(X^{(a)})^2$ of 
$d_z\theta(z-w+\delta -\Delta)$ (which is a $1$-form 
in $z$ and a function in $w$). 

This $1$-form is proportional to the lift to $X^{(a)}$ 
of $\omega$: we have 
$$
\alpha(z)dz = \kappa \omega, 
$$ 
with $\kappa \neq 0$. 
\end{lemma}

{\em Proof.}
Let us study the transformation properties of
$\alpha(z) dz$ when $z$ is transformed to 
$\gamma_{b_i}(z)$. We have 
$$
\theta(\gamma_{b_i}(z) - \gamma_{b_i}(w) + \delta - \Delta)
= 
e^{\int_z^w \omega_i}\theta(z-w + \delta - \Delta) , 
$$
therefore 
$$
d_{z} \theta(\gamma_{b_i}(z) - \gamma_{b_i}(w) + \delta - \Delta)
= 
e^{\int_z^w \omega_i} 
d_z\theta(z-w + \delta - \Delta) 
+
d_z (e^{\int_z^w \omega_i})\theta(z-w + \delta - \Delta) ;  
$$
since $\theta(z-w + \delta - \Delta)$ vanishes for
$z=w$ and $e^{\int_z^w \omega_i}$ is equal to $1$
for $z=w$, we obtain
$$
d_{z} \theta(\gamma_{b_i}(z) - \gamma_{b_i}(w) + \delta - \Delta)_{|z=w}
= 
d_z\theta(z-w + \delta - \Delta)_{|z=w} , 
$$
so that $\alpha(z)dz$ is invariant under all $\gamma_{b_i}$, and is
therefore the lift of some $1$-form $\wt \alpha(z)dz$ defined on $X$. 

Let us now determine this $1$-form. $\wt \alpha(z)dz$ is obviously 
regular on $X$. On the other hand, we have the following 
expansion of $\theta(z-w+\delta-\Delta)$ for $z$ and $w$ at the 
vicinity of some $\delta_i$ (see \cite{Mum,Fay}): 
$$  
\theta(z-w+\delta-\Delta) = z_i w_i (z_i - w_i) a(z_i,w_i), 
$$
where $z_i,w_i$ are the coordinates of $z$ and $w$ 
at $\delta_i$ and $a(z,w)$ is regular and non-zero at $(0,0)$. 
We then have 
$$  
d_z \theta(z-w+\delta-\Delta) = z_i w_i a(z_i,w_i) dz_i +  
(z_i - w_i)  \left( z_i w_i 
a_{z_i}(z_i,w_i) + w_i a(z_i,w_i) \right) , 
$$
so that 
$$  
d_z \theta(z-w+\delta-\Delta)_{| z=w} = z_i^2 a(z_i,z_i) 
dz_i , 
$$
and $\alpha(z)dz$ has a double pole at $\delta_i$. It follows that 
$\wt \alpha(z)dz$ also has a double pole at $\delta_i$, and is 
therefore proportional to $\omega$. 
\hfill \qed \medskip 

Since the function $(z,w)\mapsto \theta(z-w+ \delta - \Delta)$
vanishes on the diagonal $z=w$, it follows from this Lemma that 
$\theta(z-w+ \delta - \Delta)$ is equivalent to 
$\kappa \int_z^w \omega$ near $z=w$. 
When $\pa_h \theta(\delta - \Delta)$ is not equal to 
zero, $G_h(z,w)$ is then equivalent to 
$\pa_h \theta(\delta - \Delta) / (\kappa\int_z^w \omega)$,
whence (\ref{expansion}), with 
$C(h) = \kappa^{-1}\pa_h \theta(\delta - \Delta)$. 

Before we study the vanishing of $C(h)$, we show the following
lemma: 

\begin{lemma} \label{independence}
The $V_{\delta_i}, i = 1,\cdots, g-1$, are independent 
vectors of $\CC^g$; for $x$ a generic point of $X$, 
$V_x$ does not belong to 
$\oplus_{i=1}^{g-1} \CC V_{\delta_i}$.  
\end{lemma}

{\em Proof.}
Consider the map $\sigma : X^{(g)} \to \Jac^g(X)$, defined by 
$\sigma((y_i)) = \sum_i y_i$. By \cite{Fay}, p. 6, and
\cite{May,Lew}, this map has rank $g$ at the point $(y_i)$ iff
$h^1(\sum_i y_i) = 0$. This is the case for the point 
$\sum_i \delta_i + x$, for $x$ a generic point of $X$.  
Indeed, $h^1(\sum_i \delta_i + x)$ is then equal, by Serre 
duality, to $h^0(\sum_i \delta_i - x)$, which is zero for $x$
generic (because we have $h^0(\sum_i \delta_i)$ is equal to $1$). 
The tangent space to the image of $\sigma$ at this point is the 
span of $V_x$ and the $V_{\delta_i}$. It follows that these 
vectors are independent, for $x$ generic. 
\hfill \qed \medskip

Let us now study the vanishing of $C(h)$. 
$C(h)$ vanishes for $h$ equal 
to some $V_{\delta_i}$, because we have 
$\pa_{V_{\delta_i}} \theta(\delta - \Delta) = 
(d/dt)\theta(\sum_{j\neq i}\delta_j + \delta_i(t) 
- \Delta)_{|t=0}$, 
where $t\mapsto \delta_i(t)$ is some coordinate map from the 
vicinity of $0$ to that of $\delta_i$; on the other hand, 
$\theta(\sum_{j\neq i}\delta_j + \delta_i(t) 
- \Delta)$ vanishes identically, because of (\ref{Riem}). 

On the other hand, $C(h)$ does not vanish for $h = V_Q$, $Q$ some
point of $X$ distinct from the $\delta_i$. Indeed, 
we have 
$\pa_{V_{Q}} \theta(\delta - \Delta) = 
(d/dt)\theta( Q(t) - Q + 
\delta - \Delta)_{|t=0}$, where $t\mapsto Q(t)$ is
some coordinate map from the 
vicinity of $0$ to that of $Q$; from Lemma \ref{diagonal}
now follows that $\pa_{V_{Q}} \theta(\delta - \Delta)$ is
equal to $\omega_Q/ dt$ and is therefore not zero. 

In view of the first part of Lemma \ref{independence}, 
it follows that the linear form $C(h)$ vanishes iff  
$h$ belongs to $\oplus_{i=1}^{g-1} \CC V_{\delta_i}$, whence 
the last part of the proposition.  

\hfill \qed \medskip 

\begin{remark} \label{rem:adeles}
It follows from the proof above that 
the second statement of Lemma \ref{independence}
can be precised: for $Q$ a point of $X$, distinct from 
the $\delta_i, i = 1,\cdots, g-1$,  
$V_Q$ does not belong to 
$\oplus_{i=1}^{g-1} \CC V_{\delta_i}$.  

This statement can be translated as follows: 
let for any point $x$ of $X$, $z_x$ be some local 
coordinate at $x$. The adeles ring $\AAA$ of $X$ is
the restricted product of the formal series fields
$\CC((z_x))$. The function field $\CC(X)$ of $X$ is
embedded in $\AAA$ by taking Laurent expansions of 
a function at each point of $X$. We denote by 
$\cO_{\AAA}$ the subring of $\AAA$ of integral adeles, 
from by the restricted product of the formal series 
rings $\CC[[z_x]]$. the first cohomology ring 
$H^1(X,\cO_X)$ is defined by $H^1(X,\cO_X) 
= \AAA / \CC(X) + \cO_{\AAA}$; it is a vector space
of dimension $g$. 
For $x$ in $X$, let us denote by $z_x^{-1}$ the
element of $\AAA$ with $x$-component $z_x^{-1}$ 
and other components zero. 

Then the classes of the elements $z_{\delta_i}^{-1}$ 
and $z_x^{-1}$ form a basis of 
$H^1(X,\cO_X)$.
\hfill \qed \medskip 
\end{remark}

Our aim is now to first prove that for 
$h$ not in the span of the $V_{\delta_i}$, $L_{\CC h}$
and $\cO$ are supplementary, and then that 
$G_h$ is the corresponding Green function.  

Let us expand $G_h(z,w)$ for $w$ near $S'$. We obtain some 
element $G_h$ of $\cO \otimes (L_{\CC h} +\CC 1)$. 
This element satisfies  
$$
G_h + G_h^{(21)} = \kappa C(h)
\sum_i \eps^i \otimes \eps_i ,   
$$
for $\eps^i,\eps_i$ dual bases of $\cK$ associated with $\omega$; 
indeed, we have 
$G_h + G_h^{(21)} = C(h)
[(\int_z^w \omega)^{-1} + (\int_w^z \omega)^{-1}] = 
C(h)\kappa \al(z)^{-1}\delta(z-w)
$; on the other hand, recall that $\omega = \kappa^{-1}\alpha(z)dz$, 
so that $\sum_i \eps^i \otimes \eps_i = \al(z)^{-1}\delta(z-w)$. 
  
This implies that $L_{\CC h} + \cO = \cK$. Now, since both 
$L_{\CC h}$ and $\cO$ are isotropic and since the scalar 
product on $\cK$ is non-degenerate, their 
intersection is reduced to zero; therefore 
we have shown that $L_{\CC h}$ and $\cO$ are 
supplementary. 

Let us denote by $\bar G_h$ the Green function associated 
with this decomposition. $G_h$ is an element of $\cO \otimes 
L_{\CC h}$, and it satisfies 
$$
\bar G_h + \bar G_h^{(21)} = \sum_i \eps^i \otimes \eps_i . 
$$ 
Then the difference between $G_h - \kappa C(h)\bar G_h$ is
antisymmetric, and it belongs to $(L_{\CC h}+\CC 1)\otimes\cO$.
Since the intersection of $L_{\CC h}+\CC 1$ and $\cO$ is 
reduced to the constants, this difference is equal to zero.

Therefore: 

\begin{thm} \label{green} For $h\notin \oplus_{i} \CC V_{\delta_i}$, 
$L_{\CC h}$ and $\cO$ are supplementary. 
The Green function associated with the 
Lagrangian decomposition $\cK = \cO \oplus L_{\CC h}$
is 
$$
\bar G_h(z,w) = (\kappa C(h))^{-1}
(\pa_h\theta/\theta)(z-w+ \delta - \Delta). 
$$
\end{thm}

\begin{remark} In the case when $V = 0$, $L_V$ consists of 
the rational functions on $X$, regular outside $S'\cup \delta$ 
and with simple poles at the $\delta_i$. 
We then have $L_0 \cap \cO = \CC 1$. Indeed, this intersection
consists of the rational functions on $X$ with at most simple 
poles at $\delta$. This space is exactly $H^0(X,\cL_{\delta})$, 
which is $1$-dimensional, and therefore consists of the 
constants. 
\end{remark}

\subsection{Quantum algebras at level zero.} 
\label{level:0:theta:rels}

\subsubsection{Quasi-Hopf algebra $U_{\hbar,h}\G_{S'}$}

Let for any $s$ of $S'$, $z_s$ be a local coordinate at $s$
and $\pa_{z_s}$ be the derivation $d/dz_s$. In what follows, 
we will denote by $(z_s)$ the point of $X$ with coordinate
$z_s$. 

Let us set, for $s$ in $S'$, $i\ge 0$, 
$l_i^{(s)}(w) = {1\over {i!}} 
\pa_{z_s}^i(\pa_h\theta/\theta)
(z-w + \delta - \Delta)_{z=s}$. Then Thm. \ref{green} 
implies that $(l_i^{(s)})_{i\ge 0,s\in S'}$ is a basis of 
$L_{\CC h}$ dual to the basis $(z_s^i)_{i\ge 0,s\in S'}$ 
of $\cO_S'$. 

\begin{prop} \label{G_S:0}
Let $U_{\hbar,h}\G_{S'}$ be the algebra 
with generators $x[z_s^i],x[l_i^{(s)}]$, $s\in S', i\ge 0$, 
$x = e,f,h$, generating series 
$$
x^{(s)}(z_s) = \sum_{i\ge 0} x[l_i^{(s)}] z_s^i + 
\sum_{i\ge 0,t\in S'} x[z_t^i] l_i^{(t)}((z_s)), \quad 
x = e,f,
$$ 
$$
h^+(z) = \sum_{i\ge 0} h[z_s^i] l_i^{(s)}(z), 
\quad 
h^{-(s)}(z_s) = \sum_{i\ge 0} h[l_i^{(s)}] z_s^i , 
$$
and relations
\begin{equation} \label{h:h:level0:g}
[h[\alpha] , h[\beta]] = 0, 
\end{equation}
for any $\al,\beta$, 
\begin{equation} \label{K+:e:level0:g}
K^+(z) e^{(s)}(w_s) K^+(z)^{-1} = 
{{\theta(z-(w_s) - \hbar h + \delta - \Delta)} 
\over {\theta(z-(w_s) + \hbar h + \delta - \Delta)}} 
e^{(s)}(w_s), 
\end{equation}
\begin{equation} \label{K-:e:level0:g}
K^{-(s)}(z_s) e^{(t)}(w_t) K^{-(s)}(z_s)^{-1} = 
{{\theta((z_s)-(w_t) + \hbar h + \delta - \Delta)} 
\over {\theta((z_s)-(w_t) - \hbar h + \delta - \Delta)}}  
e^{(t)}(w_t), 
\end{equation}
\begin{equation} \label{K+:f:level0:g}
K^+(z) f^{(s)}(w_s) K^+(z)^{-1} = {{\theta(z-(w_s) + \hbar h 
+ \delta - \Delta)} 
\over {\theta(z-(w_s) - \hbar h + \delta - \Delta)}}  
f^{(s)}(w_s), 
\end{equation}
\begin{equation} \label{K-:f:level0:g}
K^{-(s)}(z_s) f^{(t)}(w_t) K^{-(s)}(z_s)^{-1} = 
{{\theta((z_s)-(w_t) - \hbar h + \delta - \Delta)} 
\over {\theta((z_s)-(w_t) + \hbar h + \delta - \Delta)}} f^{(t)}(w_t), 
\end{equation}
\begin{align} \label{e:e:level0:g}
& 
\theta((z_s)-(w_t)+ \hbar h + \delta - \Delta) 
e^{(s)}(z_s)e^{(t)}(w_t)
\\ & \nonumber 
=
\theta((z_s)-(w_t) - \hbar h + \delta - \Delta) 
e^{(t)}(w_t) e^{(s)}(z_s), 
\end{align}
\begin{align} \label{f:f:level0:g}
& 
\theta((z_s)-(w_t) - \hbar h + \delta - \Delta) 
f^{(s)}(z_s)f^{(t)}(w_t)
\\ & \nonumber 
=
\theta((z_s)-(w_t) +\hbar h + \delta - \Delta) 
f^{(t)}(w_t) f^{(s)}(z_s), 
\end{align}
and
\begin{equation} \label{e:f:level0:g}
[e^{(s)}(z_s) , f^{(t)}(w_t)] = \delta_{st}
\delta(z_s-w_t) \left( K^+((z_s)) - K^{-(s)}(z_s)^{-1} \right) , 
\end{equation}
with $K^{+}(z), K^{-(s)}(z_s)$ defined as in sect. \ref{level:0},  
coproduct $\Delta_{S',\delta}$ defined by (\ref{copdt:1}), (\ref{copdt:2}), 
is a quantization of the double Lie bialgebra
structure on $\bar\G\otimes \cK_{S'}$ defined by the decomposition 
$$
\bar\G\otimes \cK_{S'} 
= (\bar\HH\otimes \cO_S \oplus \bar\N_+\otimes\cK_{S'}) 
\oplus 
(\bar\HH\otimes L_{\CC h} \oplus \bar\N_-\otimes\cK_{S'}).
$$
Thm. \ref{level:0:qH} can be applied to it to define a 
quantization 
of the double quasi-Lie bialgebra structure on $\bar\G\otimes \cK_{S'}$
defined by the decomposition
$$
\bar\G\otimes \cK_{S'}
= (\bar\G\otimes \cO_{S'} ) \oplus (\bar\G\otimes L_{\CC h}).
$$
\end{prop}

{\em Proof.} This follows from Thm. \ref{green}, the expansion 
$$
{{\theta(z-w+\hbar h + \delta - \Delta)}\over
{\theta(z-w + \delta - \Delta)}} 
= a(z,w) + \hbar {{\pa_h\theta}\over{\theta}}(z-w + \delta - \Delta) b(z,w), 
$$
where $a(z,w)$ and $b(z,w)$ belong to $\cO\hat\otimes\cO[[\hbar]]$, 
and Thm. \ref{level:0:qH}. 
\hfill \qed \medskip

\subsubsection{Algebra $U_{\hbar,h}\G_{S',\delta}$}

Let $L_{\CC h}^{S',\delta}$ be the space of functions $f$ 
defined on 
$\wt X$, regular outside the lifts of $S'$ and $\delta$ 
(here and later, we will also denote by $\delta$ the 
support $\{\delta_i\}$ of $\delta$), 
such that the differences $f(\gamma_{b_i}z) - f(z)$ are constant
and form a vector proportional to $h = (h_i)_{1\leq i\leq g}$, 
and such that 
$$
\sum_{i=1}^g h_i \left( \int_{a_i}f\omega + 
\int_{\gamma_{b_i}(a_i)}f\omega \right) = 0. 
$$
On the other hand, let $\wt\cO_{S',\delta}$ be the direct sum 
$\cO_{S'} \oplus 
(\oplus_{i=1}^{g-1} z_{\delta_i}^{-1} \cO_{\delta_i})$. 

\begin{prop} \label{what:is:green}
Endow $\cK_{S',\delta}$ with the scalar product defined by 
$\langle \phi,\psi \rangle_{\cK_{S',\delta}} 
= \sum_{\al\in S'\cup\delta} $ 
$\res_{\al}(\phi\psi\omega).$ 
The spaces $L_{\CC h}^{S',\delta}$ and 
$\wt\cO_{S',\delta}$ are isotropic supplementaries in 
$\cK_{S',\delta}$. The Green function associated to this 
decomposition is given by the collection of expansions, for
$w$ near each point of $S'\cup\delta$, of 
$$
\wt G_h(z,w) = {{\pa_h\theta}\over{\theta}}(z-w + \delta - \Delta). 
$$ 
\end{prop}

{\em Proof.} The argument showing that $L_{\CC h}^{S',\delta}$
is isotropic is similar to the argument used for   
$L_{\CC h}$. On the other hand, since $\omega$ is regular 
on $S'$ and has double poles at the $\delta_i$, 
$\wt\cO_{S',\delta}$ is also isotropic. From Rem. \ref{rem:adeles}
also follows that the direct sum of these spaces is 
$\wt\cO_{S',\delta}$. This proves the first part of the 
proposition. 

To prove its second part, let us expand $\wt G_h(z,w)$ for
$w$ near each point of $S'\cup \delta$. Since for fixed $z$, 
the function $w\mapsto \theta(z - w + \delta - \Delta)$
either vanishes to first order (for $w$ near $\delta_i$)
or is non-zero (for $w$ near $S'$), this expansion will 
be a series $\sum_{_la} f_{\la} \otimes o_{\la}$, with 
$o_{\la}$ in $\wt\cO_{S',\delta}$. On the other hand, as a
function of $z$,    
$\wt G_h(z,w)$ is regular for $z$ outside $\delta$ and $w$, 
and has the functional properties (\ref{q-per}), so that 
the $f_{\la}$ belong to $L_{\CC h}^{S',\delta}\oplus \CC 1$. 

Let us compare now the resulting expansion of $\wt G_h(z,w)$ 
with the Green function $G_{S',\delta}$ of the decomposition 
$L_{\CC h}^{S',\delta} \oplus \wt\cO_{S',\delta}$. 
As is Prop. \ref{what:is:green}, we can check that the sums 
$\wt G_h + \wt G_h^{(21)}$ and $G_{S',\delta} + 
G_{S',\delta}^{(21)}$ coincide with the same delta-functions. 
We conclude from there that the difference 
$\wt G_h - G_{S',\delta}$ is antisymmetric. Since it belongs
to the tensor square of the intersection of 
$\wt\cO_{S',\delta}$ and 
$L_{\CC h}^{S',\delta}\oplus \CC 1$, which is $\CC 1$, 
this difference is zero. 
\hfill \qed\medskip 

Let us define now $(U_{\hbar,h}\G_{S',\delta}, \Delta)$ 
as the algebra defined by the generators and relations of 
Prop. \ref{G_S:0}, with $S'$ replaced by $S'\cup\delta$, 
$L_{\CC h},\cO_{S'}$ by $L_{\CC h}^{S',\delta}$ and 
$\wt\cO_{S',\delta}$. 

\begin{lemma}
$U_{\hbar,h}\G_{S',\delta}$ 
is a flat deformation of the 
enveloping algebra of $\G_{S',\delta} = \bar\G\otimes 
\cK_{S',\delta}$. 
\end{lemma}

{\em Proof.} We first prove: 
\begin{lemma} The subalgebra of $U_{\hbar,h}\G_{S',\delta}$  
generated by the $e[\phi],\phi\in \cK_{S',\delta}$
is a flat deformation of the corresponding classical 
subalgebra. 
\end{lemma}

{\em Proof.} Let us first consider the subalgebras $A_{\al}$
generated by the $e[\phi],\phi\in \cK_{\al}, \al\in 
S'\cup\delta$. 
The function $(z,w)\mapsto \theta(z-w + \delta - \Delta)$ has 
the following behavior:  
for $z,w$ near some $\delta_i$, we have 
$$
\theta(z-w + \delta - \Delta) = z_i w_i (z_i-w_i) \phi_i(z_i,w_i), 
$$
with $\phi_i(z_i,w_i)$ invertible in $\CC[[z_i,w_i]][[\hbar]]$; 
for $z,w$ near some $s\in S'$, we have
$$
\theta(z-w + \delta - \Delta) = (z_s-w_s) \phi_s(z_s,w_s), 
$$ 
with $\phi_s(z_s,w_s)$ invertible in $\CC[[z_s,w_s]][[\hbar]]$. 
After we divide them by $z_i w_i \phi(z_i,w_i)$ in the first case, 
and by $\phi(z_s,w_s)$ in the second case, 
the relations between the fields $e^{(\al)}(z_\al)$ 
have the form of the vertex relations of Prop. \ref{PBW}. 
It follows that each algebra $A_{\al}$
is a flat deformation of the corresponding classical 
subalgebra. 

Let us now study the relations between the various 
$e[\phi_\al], \phi_\al \in \cK_{\al}$. Recall that 
the function
$(z,w)\mapsto \theta(z-w + \delta - \Delta)$ has simple 
zeroes for $w$ equal to $z$, or for $z$ or $w$ equal to 
one of the $\delta_i$. It follows that, after we divide 
the relation between $e^{(\al)}(z_{\al})$ and 
$e^{(\beta)}(z_{\beta})$, by $z_{\al}$ if $\alpha$
belongs to $\delta$, and by $z_{\beta}$ if $\beta$
belongs to $\delta$, we obtain a relation of the form 
$$
e^{(\al)}(z_\al)e^{(\beta)}(w_\beta) = 
f_{\al\beta}(z_\al,w_\beta) e^{(\beta)}(w_\beta)
e^{(\al)}(z_\al), 
$$
with $f_{\al\beta}(z_\al,w_\beta)$ invertible in 
$\CC[[z_\al,w_\beta]][[\hbar]]$. It follows that 
monomials in the $e[\phi_\al], \phi_\al \in \cK_{\al}$, can 
be expressed as sums of monomials with the $\alpha$ occurring
in a prescribed order. 
\hfill \qed \medskip 

Therefore the subalgebra of $U_{\hbar,h}\G_{S',\delta}$  
generated by the $e[\alpha],\alpha\in \cK_{S',\delta}$
is a flat deformation of the corresponding classical 
subalgebra. 
One may then obtain the analogous result 
for the subalgebra generated by the fields $f[\alpha]$. 
The result then follows from the triangular decomposition 
of $U_{\hbar,h}\G_{S',\delta}$. 
\hfill \qed \medskip

The Hopf algebra $(U_{\hbar,h}\G_{S',\delta} , \Delta)$ is then the 
quantization of the double structure on 
$\bar\G \otimes \cK_{S',\delta}$ given by the 
decomposition
$$
\bar\G \otimes \cK_{S',\delta} 
= 
(\bar\HH \otimes L_{\CC h}^{S',\delta} \oplus 
\bar \N_+ \otimes \cK_{S',\delta})
\oplus 
( \bar\HH \otimes \wt\cO_{S',\delta} 
\oplus 
\bar \N_- \otimes \cK_{S',\delta}) . 
$$

Since $\wt\cO_{S',\delta}$ is not a subring of $\cK_{S',\delta}$, 
we cannot expect find a corresponding subalgebra of 
$U_{\hbar,h}\G_{S',\delta}$. 
However, we have: 

\begin{prop} \label{subalg}
The subalgebra of $U_{\hbar,h}\G_{S',\delta}$ generated by the 
$x[\phi_i],\phi_i\in \cK_{\delta_i}$, and the 
$x[o_s],o_s\in \cO_s$, $x=e,f,h$, is a flat deformation of the 
enveloping algebra of $\bar\G\otimes (\cK_{\delta}\oplus \cO_{S'})$.  
\end{prop}

{\em Proof.} Let us first prove that the similar statement 
is true for the subalgebra $N_+$ of $U_{\hbar,h}\G_{S',\delta}$ 
generated by the 
$e[\phi_i],\phi_i\in \cK_{\delta_i}$, and the 
$e[o_s],o_s\in \cO_s$.
The commutation relations between the 
$e[\phi_\al],\phi_\al\in \cK_{\al}$, and the 
$e[\phi_\beta],\phi_\beta\in \cK_{\beta}$, for $\al\neq\beta$, 
are of the form 
$$
e[z_\al^n]e[z_\beta^m] 
= \sum_{n',m'\in \ZZ} a(\al,\beta)_{nm}^{n'm'} 
e[z_\beta^{m'}] e[z_\al^{n'}],  
$$
and the summation is on $n'\ge 0$ (resp. $m'\ge 0$) if $\al$ is in 
$\delta$ an $n\ge 0$ (resp. $\beta$ is in $\delta$ and $m\ge 0$).
It follows that a basis of $U_{\hbar,h}\G_{S',\delta}$ 
is given by the products of 
bases for its subalgebras generated by the $e[z_s^n], n\in\ZZ, 
s\in S'$, and the 
$e[z_i^n], n\ge 0$. These subalgebras are flat deformations
of their classical analogues: this is because the    
$e[z_s^n], n\in\ZZ$ are subject to vertex relations, and
because we have Yangian-type commutation 
relations between the $e[o_s],o_s\in \cO_s$. 

Let us now prove that the algebra $N_+$ is stable by the 
adjoint action of the $h[\phi_i],\phi_i\in \cK_{\delta_i}$, and 
of the $h[o_s],o_s\in \cO_s$. From the identity (\ref{K+:e:level0:g}) 
follows that  
\begin{equation} \label{modes:h:e}
[(1+V)h^+(z),e^{(\al)}(z_\al)] = 
\log {{\theta(z-(z_\al)+\hbar h + \delta - \Delta)}
\over {\theta(z-(z_\al)-\hbar h + \delta - \Delta)}}
e^{(\al)}(z_\al), 
\end{equation}
for $\al$ in $S'$ or $\delta$. 
We then check: 

1) the adjoint action of any $h[\phi],\phi\in \cO_s$ preserves 
the linear space spanned by the 
$e[\phi_i],\phi_i\in\cK_{\delta_i}$ and the 
$e[z_t^n],n\ge 0$, $t$ in $S'$. 
This 
is because (\ref{modes:h:e}) yields, for $\al$ is $S'$ or $\de$, 
\begin{align*}
& [h^+[z_s^n],e^{(\al)}(z_\al)] = 
\res_{z_s = 0} 
\left( 
\log {{\theta((z_s)-(z_\al)+\hbar h + \delta - \Delta)}
\over {\theta((z_s)-(z_\al)-\hbar h + \delta - \Delta)}}
z_s^n dz_s \right) 
e^{(\al)}(z_\al), 
\end{align*}
and the function in the right side is regular for $z_\al=0$ if
$\al$ is in $S'$. 

2) the adjoint action of any $h[\phi],\phi\in \cK_{\de_i}$ 
preserves the linear space spanned by the 
$e[\phi_j],\phi_j\in\cK_{\delta_j}$ and the 
$e[o_s],o_s\in \cO_s$. The first statement is clear. To 
show the second one, we first prove in a way analogous to 
point 1) above that the adjoint action of any 
$h[\phi],\phi\in z_i^{-1}\cO_i$ preserves the linear 
space formed of the $e[o_s],o_s\in \cO_s$. After that, 
(\ref{K-:e:level0:g}) implies that 
$$
[h^{-(i)}(z_i),e^{(s)}(z_s)] = 
\log {{\theta((z_i)-(z_s)+\hbar h + \delta - \Delta)
}\over{\theta((z_i)-(z_s)-\hbar h + \delta - \Delta)}}
e^{(s)}(z_s) . 
$$
Let then $(r_i^{(k)},r_s^{(k)})_{k\ge 0}$ be the dual basis in 
$L^{S',\delta}_{\CC h}$
to $(z_i^{k-1},z_s^k)_{k\ge 0}$. Each $r_i^{k}$ is then 
multivalued on $X$, it has  
poles of order $-k-2$ at $\delta_i$, of order $1$ at most 
as each $s$, and is regular at the other points. We have 
$$
[h[r_i^{(k)}],e^{(s)}(z_s)] = 
\res_{z_i=0} z_i^2 dz_i r_i^{(k)}(z_i)
\log {{\theta((z_i)-(z_s)+\hbar h + \delta - \Delta)
}\over{\theta((z_i)-(z_s)-\hbar h + \delta - \Delta)}}
e^{(s)}(z_s) , 
$$
and since the function in the right side is regular for 
$z_s = 0$, the adjoint action of  
$h[r_i^{(k)}]$ 
preserves $\{ e[o_s],o_s\in \cO_s \}$.
Since any element $\cK_i$ can be obtained by linear
combination of the $h[r_i^{(k)}]$ and of the $h[\la]$, 
$\la\in (\oplus_i z_{\delta_i}^{-1}\cO_{\delta_i}) 
\oplus (\oplus_{s\in S'}\cO_s)$, the same result is true
for any $h[\la],\la\in \cK_{\delta_i}$. 

After that, we prove that the subalgebra $N_+$ of
$U_{\hbar,h}\G_{S',\delta}$ generated by the $f[\phi]$, 
$\phi$ in $\cO_s$ and in $\cK_{\delta_i}$, is a flat deformation
of its classical analogue, using the same reasoning as for 
$N_+$. Finally, any commutator $[e[\phi],f[\psi]]$, $\phi,\psi$ 
in the sums of $\cO_s$ and $\cK_{\delta_i}$, is expressed as 
products of $h[\phi]$, $\phi$ in $\cO_s$ or $\cK_{\delta_i}$. 
\hfill \qed \medskip  

\begin{remark} Following the proof of Thm. \ref{isomorphisms}, 
one can show that the algebra 
$U_{\hbar,h}\G_{S',\delta}$ is isomorphic to a quotient 
(with central elements identified) of the tensor product 
$DY(\SL_2)_0^{\prime\otimes \card S'} \otimes 
U_{\hbar,h}\G_{S'}$,  where $U_{\hbar,h}\G_{S'}$ is the
algebra corresponding to an empty $S'$, and $DY(\SL_2)_0'$ 
is the double Yangian 
algebra (without derivation nor central element). The 
subalgebra of Prop. \ref{subalg} could then be identified with 
$Y(\SL_2)^{\otimes \card S'} \otimes 
U_{\hbar,h}\G_{S'}$, where $Y(\SL_2)$ is the Yangian subalgebra of 
$DY(\SL_2)'_0$ generated by the nonnegative modes generators.
\end{remark}

\subsubsection{Regular subalgebra in $U_{\hbar,h}\G_{S',\delta}$}
\label{regular:level0:sect}

Let us define, for $\eps$ in $\cK_{S',\delta}$, $\bar h[\eps]$
as the Cartan element of $U_{\hbar,h}\G_{S',\delta}$ such that 
$$
[\bar h[\eps],e^{\al}(z_{\al})] = 2\eps^{(\al)}(z_{\al}) 
e^{\al}(z_{\al}), 
$$
for any $\al$ in $S'\cup\delta$. 
For $\phi$ a regular function on 
$\wt X - \wt S'$, let us set $\bar h[\phi] = \sum_{\al\in S'\cup\delta} 
\bar h[\phi^{(\al)}]$, where $\phi^{(\al)}$ is the image in $\cK_{S',\delta}$
of the element of $\cK_{\al}$ given by 
expansion of 
$\phi$ near $\wt \al$

Let $P$ be a point of $S'$. 
Let set  
$$
x(z) = \exp \left( {1\over 4}\log{{\theta(z- P+\hbar h  
+ \delta - \Delta)}
\over{\theta(z - P -\hbar h + \Delta + \delta)}} 
\bar h[1] - {1\over 4} \bar h[\log{{\theta(\cdot - P
+\hbar h + \delta - \Delta)}
\over{\theta(\cdot - P -\hbar h + \delta - \Delta)}} ]\right),
$$
and 
\begin{equation} \label{bar:1}
\bar e^{(\al)}(z_\al) = e^{(\al)}(z_\al) x(z_\al), 
\quad
\bar f^{(\al)}(z_\al) = f^{(\al)}(z_\al) x(z_\al), 
\end{equation}
\begin{equation} \label{bar:2}
\bar K^{+}(z) = K^+(z)x(z)^2 , 
\quad 
\bar K^{-(\al)}(z_{\al}) 
= K^{-(\al)}(z_{\al}) x((z_{\al}))^2. 
\end{equation} 

\begin{lemma} \label{rels:reg}
Let us set
$$
q_0(z,w) = {{\theta(z - w + \hbar h + \delta - \Delta)} 
\over 
{\theta(z - w - \hbar h + \delta - \Delta)}}, \quad
q_{\pm}(z,w) = \theta(z - w \pm \hbar h + \delta - \Delta), 
$$
and
$$
q(z,w) = {{q_0(z,w)}\over{q_0(z,P)q_0(P,w)}}, 
\quad
\wt q_{\pm}(z,w) = 
{{q_{\pm}(z,w)}\over{q_{\pm}(z,P)q_{\pm}(P,w)}}; 
$$
then
$\bar e^{(\al)}(z_{\al}),\bar f^{(\al)}(z_{\al}),\bar K^{+}(z)$
and $\bar K^{-(\al)}(z_{\al})$ satisfy the relations
$$
\bar K^+(z) \bar e^{(\al)}(w_{\al}) 
\bar K^+(z)^{-1} = 
q(z,(w_{\al}))
\bar e^{(\al)}(w_{\al}) 
$$
$$
\bar K^+(z)\bar f^{(\al)}(w_{\al}) \bar K^+(z)^{-1} = 
q(z,(w_{\al}))^{-1}
\bar f^{(\al)}(w_{\al}) 
$$
$$
\bar K^{-(\al)}(z_{\al}) \bar e^{(\beta)}(w_{\beta}) 
\bar K^{-(\beta)}(z_{\beta})^{-1} = 
q((z_{\al}),(w_{\beta}))^{-1} 
\bar e^{(\beta)}(w_{\beta}) 
$$
$$
\bar K^{-(\al)}(z_{\al}) \bar f^{(\beta)}(w_{\beta}) 
\bar K^{-(\beta)}(z_{\beta})^{-1} = 
q((z_{\al}),(w_{\beta})) \bar f^{(\beta)}(w_{\beta})   
$$
$$
[\bar e^{(\al)}(z_{\al}) , \bar f^{(\beta)}(w_{\beta})] = 
\delta_{\al\beta} \delta(z_{\al} - w_{\beta}) 
\left( \bar K^+((z_{\al})) - \bar K^{-(\al)}(z_{\al}) \right), 
$$
$$
\wt q_+((z_\al),(w_\beta))
\bar e^{(\al)}(z_\al)\bar e^{(\beta)}(w_\beta)
 =
\wt q_-((z_\al),(w_\beta))
\bar e^{(\beta)}(w_\beta) \bar e^{(\al)}(z_\al), 
$$
$$
\wt q_-((z_\al),(w_\beta))
\bar f^{(\al)}(z_\al)\bar f^{(\beta)}(w_\beta)
=
\wt q_+((z_\al),(w_\beta))
\bar f^{(\beta)}(w_\beta) \bar f^{(\al)}(z_\al). 
$$
\end{lemma}

{\em Proof.}
This follows from the identities  
\begin{align*}
& x(z) e^{(\al)}(w_{\al}) x(z)^{-1} = \\ &
\left( 
{{\theta(z-P+\hbar h + \delta - \Delta)}
\over{\theta(z-P - \hbar h + \delta - \Delta)}} 
: 
{{\theta((w_{\al}) -P+\hbar h + \delta - \Delta)}
\over{\theta((w_{\al}) -P - \hbar h + \delta - \Delta)}} 
\right)^{1/2}
e^{(\al)}(w_{\al}) , 
\end{align*}
\begin{align*}
& x(z) f^{(\al)}(w_{\al}) x(z)^{-1} = \\ &  
\left( 
{{\theta(z -P-\hbar h + \delta - \Delta)}
\over{\theta(z-P + \hbar h + \delta - \Delta)}} 
: 
{{\theta((w_{\al}) -P-\hbar h + \delta - \Delta)}
\over{\theta((w_{\al}) -P + \hbar h  + \delta - \Delta)}} 
\right)^{1/2}
f^{(\al)}(w_{\al}) . 
\end{align*}
\hfill \qed \medskip

Let $R_{S',\delta}$ be the algebra of functions on $X$, regular
outside $S'\cup\delta$. Then the intersection of 
$R_{S',\delta}$ and $L_{\CC h}^{S',\delta}$ has codimension
$1$ in each of these spaces. A supplementary of this 
intersection in $L_{\CC h}^{S',\delta}$ is spanned by 
$e_0(z) = (\pa_h\theta/\theta)(z - P + \delta - \Delta) + c$, 
with $c$ a certain constant. 

Let us define $\Sigma$ as the direct sum of $\CC e_0$ and 
the orthogonal of $e_0$ in $\wt\cO_{S',\delta}$. Then: 

\begin{lemma}
The spaces $R_{S',\delta}$ and $\Sigma$ are isotropic 
supplementaries in $\cK_{S',\delta}$. The corresponding Green
kernel is proportional to the collection of expansions, for $w$ near 
the points of $S'\cup\delta$, of  
$$
G_R(z,w) = 
{{\pa_h\theta}\over{\theta}}(z-w+\delta-\Delta)
- 
{{\pa_h\theta}\over{\theta}}(z-P+\delta-\Delta)
+
{{\pa_h\theta}\over{\theta}}(w-P+\delta-\Delta) . 
$$
\end{lemma} 

{\em Proof.} Let us set $e_{-1}=1$; then the pairing 
$\langle e_0,e_{-1} \rangle_{\cK_{S',\delta}}$ is nonzero (otherwise 
$R_{S',\delta}+L_{\CC h}^{S',\delta}$ would be isotropic). 
Let us complete $e_0$ and $e_{-1}$ to dual bases 
$(e_i)_{i\ge 0}$ and $(e_{-i-1})_{i\ge 0}$ of $R_{S',\delta}$
and $L_{\CC h}^{S',\delta}$. Then bases for $R_{S',\delta}$ and
$\Sigma$ are $(e_{-1},e_{i})_{i>0}$ and $(e_0,e_{-i-1})_{i> 0}$, 
so that both spaces are supplementary. The difference between 
the Green function for this decomposition and $\cK_{S',\delta}
= L_{\CC h}^{S',\delta} \oplus \wt\cO_{S',\delta}$ is just 
$e_0 \otimes e_{-1} - e_{-1} \otimes e_0$. 
\hfill \qed \medskip

We can therefore construct an other double quasi-Lie 
bialgebra
structure on $\bar\G \otimes \cK_{S',\delta}$, based on 
this decomposition: 
$$
\bar\G \otimes \cK_{S',\delta}
= (\bar\G \otimes R_{S',\delta}) \oplus 
(\bar\G \otimes \wt\cO_{S',\delta})
$$
and its usual infinite twist
\begin{equation} \label{double:reg}
\bar\G \otimes \cK_{S',\delta}
= (\bar\HH \otimes R_{S',\delta}
\oplus \bar\N_+ \otimes \cK_{S',\delta}
) \oplus 
(\bar\HH \otimes \wt\cO_{S',\delta}
\oplus \bar\N_- \otimes \cK_{S',\delta}) . 
\end{equation}

Set 
$q_{0+}(z,w) = {{\theta(z - w -\hbar h + \delta - \Delta)}
\over{\theta(z - w + \delta - \Delta)}}$. We have then 
\begin{equation} \label{note1}
q_{0+}(\gamma_{a_i}z,w) = q_{0+}(z,\gamma_{a_i}w) = q_{0+}(z,w), 
\end{equation}
\begin{equation} \label{note2}
q_{0+}(\gamma_{b_i}z,w) = e^{-h_i}q_{0+}(z,w), \quad 
q_{0+}(z,\gamma_{b_i}w) = e^{h_i}q_{0+}(z,w), 
\end{equation}
so that the function $(z,w)\mapsto {{q_{0+}(z,w)}\over
{q_{0+}(z,P)q_{0+}(P,w)}}$ is single-valued on the complement
of the diagonal of $X - (P\cup\delta)$. 

We then have: 

\begin{lemma} \label{functions}
Set 
$q_{0+}(z,w) = {{\theta(z - w -\hbar h + \delta - \Delta)}
\over{\theta(z - w + \delta - \Delta)}}$. Then 
for some $a,b$ in
$R_{S',\delta}^{\otimes 2}[[\hbar]]$, we have 
$$
{{q_{0+}(z,w)}\over{q_{0+}(z,P) q_{0+}(P,w)}} = 
a(z,w) + b(z,w) G_R(z,w). 
$$
\end{lemma}

{\em Proof.} For $z$ close to $w$, we have the expansions
$G_P(z,w) = C(h) / \int_z^w \omega + O(1)$, and 
$$
{{q_{0+}(z,w)}\over {q_{0+}(z,P)q_{0+}(P,w)}} = 
{1\over{q_{0+}(z,P)q_{0+}(P,z)}} {{\theta(\hbar h)}
\over{\kappa \int_z^w \omega}} + O(1), 
$$
by the remark following Lemma \ref{diagonal}. 

Set then 
$$
b(z,w) = {{\theta(\hbar h) / (C(h)\kappa)}\over{q_{0+}(z,P)
q_{0+}(P,z)}} ; 
$$
by (\ref{note1}) and (\ref{note2}), this function is single-valued
on $X$; as its only poles are at $P$ and $\delta$,  
this is a series in $\hbar$ with coefficients in $R_{S',\delta}^{\otimes 2}$. 
Set then
$$
a(z,w) = {{q_{0+}(z,w)}\over {q_{0+}(z,P)q_{0+}(P,w)}} 
- b(z,w) G_R(z,w) ; 
$$ 
this is again a single-valued function on the complement of the 
diagonal of $(X - (\delta \cup P))^2$, which is also regular on the 
diagonal. 
\hfill \qed \medskip

According to Thm. \ref{level:0:qH}, we can then   
define a quantization 
$(U_{\hbar,h}\wt\G_{S',\delta},\wt\Delta_{S',\delta})$ 
of the quasi-Lie bialgebra structure on 
$\bar\G \otimes \cK_{S',\delta}$ 
defined by the decomposition (\ref{double:reg}), 
using the functions $a$ and $b$ of Lemma \ref{functions}. 
Denote by $\wt e^{(\al)}(z_{\al}),\wt f^{(\al)}(z_{\al}),
\wt K^{+}(z)$ and $\wt K^{-(\al)}(z_{\al})$ the 
generating fields of this algebra, analogues to the fields 
$e(z),f(z)$, $K^{-}(z)$ and $K^{+}(z)$ of sect. \ref{level:0}
(note the inversion of indices of fields $K$).

\begin{prop} \label{isom:i}
The map $i$ assigning to the fields 
$\bar e^{(\al)}(z_{\al}),\bar f^{(\al)}(z_{\al}),
\bar K^{+}(z)$ and $\bar K^{-(\al)}(z_{\al})$ defined
by (\ref{bar:1}) and (\ref{bar:2}), the fields
$\wt e^{(\al)}(z_{\al}),\wt f^{(\al)}(z_{\al}),
\wt K^{+}(z)$ and $\wt K^{-(\al)}(z_{\al})$ respectively, 
defines an algebra isomorphism from $U_{\hbar,h}\G_{S',\delta}$ 
to $U_{\hbar,h}\wt\G_{S',\delta}$. 
The coproducts $\Delta_{S',\delta}$ and $\wt\Delta_{S',\delta}$
are then connected by a twist transformation
$$
\wt\Delta_{S',\delta}(i(x)) = F_0 (i\otimes i)(\Delta_{S',\delta}(x))
F_0^{-1}, 
$$
for any $x$ in $U_{\hbar,h}\G_{S',\delta}$, where 
\begin{align*}
F_0 & = \exp [{1\over 8}(\bar h[\log 
{{\theta(\cdot - P +\hbar h +\delta - \Delta)}
\over {\theta(\cdot - P -\hbar h +\delta - \Delta)}}] 
\otimes \bar h[1] \\ & 
- \bar h[1] 
\otimes \bar h[\log {{\theta(\cdot - P +\hbar h +\delta - \Delta)}
\over {\theta(\cdot - P -\hbar h +\delta - \Delta)}}] )].  
\end{align*}
\end{prop}

{\em Proof.}
Let us first check that the relations defining $i$ are consistent. 
For $\eps$ in $\cK_{S',\delta}$, let $\wt h[\eps]$ be the Cartan element
of $U_{\hbar,h}\wt\G_{S',\delta}$, such that 
$$
[\wt h[\eps],\wt e^{\al}(z_{\al})] = 
\eps^{(\al)}(z_{\al}) 
\wt e^{\al}(z_{\al}). 
$$
For some functions $\la(z,w_\beta), \la^{(\al)}(z_{\al},w_\beta)$, 
we have 
$$
[\log \bar K^+(z),e^{(\beta)}(w_{\beta})] = \la(z,w_\beta) 
e^{(\beta)}(w_{\beta})
,
$$
$$ 
[\log \bar K^{-(\al)}(z_{\al}),e^{(\beta)}(w_{\beta})] = 
\la^{(\al)}(z_{\al},w_\beta) e^{(\beta)}(w_{\beta}) , 
$$
as well as 
$$
[\log \wt K^+(z),\wt e^{(\beta)}(w_{\beta})] = \la(z,w_\beta) 
\wt e^{(\beta)}(w_{\beta})
,
$$
$$ 
[\log \wt K^{-(\al)}(z_{\al})
,\wt e^{(\beta)}(w_{\beta})] = 
\la^{(\al)}(z_{\al},w_\beta) \wt e^{(\beta)}(w_{\beta}) . 
$$
Expand $\la(z,w_\beta) = \sum_i a_i(z) b_i(w_\beta)
, \la^{(\al)}(z_{\al},w_\beta) = \sum_i a^{(\al)}_i(z_\al) 
b^{(\al)}_i(w_\beta)$;  
we have then 
$$
\log \bar K^+(z) = \sum_i  \bar h[b_i] a_i(z) , \quad 
\log \bar K^{-(\al)}(z_{\al}) = \sum_i \bar h[b_i^{(\al)}] 
a_i^{(\al)}(z_\al), 
$$
and
$$
\log \wt K^+(z) = \sum_i  \wt h[b_i] a_i(z) , \quad 
\log \wt K^{-(\al)}(z_{\al}) = \sum_i \wt h[b_i^{(\al)}] 
a_i^{(\al)}(z_\al). 
$$ 
It follows that the generating formulas
for $i(\bar K^+(z)),i(\bar K^{-(\al)}(z_{\al}))$ are consistent 
and yield $i(\bar h[\eps]) = \wt h[\eps]$, 
for any $\eps$ in $\cK_{S',\delta}$. 

The relations of Lemma \ref{rels:reg} then imply that $i$ is an 
algebra morphism. The conjugation identity is checked 
directly. \hfill \qed \medskip

Set for $x = \bar e,\bar f,\bar h$, and $\phi\in\cK_{S',\delta}$, 
$x[\phi] = \sum_{\al\in S'\cup\delta} \res_{\al}(x(z)\phi(z)\omega)$. 
Then: 

\begin{corollary} \label{lancaster}
The subalgebra of $U_{\hbar,h}\G_{S',\delta}$
generated by the $\bar e[r],\bar f[r]$ and $\bar K^+[r]$, for $r$
in $R_{S',\delta}$, is a flat deformation of the enveloping algebra
of $\bar\G\otimes R_{S',\delta}$. 
\end{corollary}

{\em Proof.} After we apply $i$, this is follows from the PBW result 
of Thm. \ref{level:0:qH} on regular subalgebras of the algebras 
$U_{a,b}\G$. \hfill \qed \medskip 

\begin{remark} The field $K^+(z)$ satisfies the functional 
equations
$$
K^+(\gamma_{a_i}z) = K^+(z), \quad
K^+(\gamma_{b_i}z) = K^+(z) \bar h[1]^{-\hbar h_i}, 
$$
analogous to relation (44) of \cite{ell:QG}. 
\end{remark} 

\begin{remark} It is easy to specialize Prop. \ref{isom:i}
to obtain an isomorphism between the centerless versions 
of the elliptic algebras of \cite{Enr-Rub} and \cite{ell:QG}.  
\end{remark}

\begin{remark} \label{g>1:nonzero:level}
By analogy with \cite{ell:QG}, one may construct a 
``centrally extended'' version $U_{\hbar,h}\hat\G_{S',\delta}$, 
of the algebra $U_{\hbar,h}\G_{S',\delta}$, 
with additional central generator 
$K$ and relations (\ref{h:h:level0:g}), (\ref{e:f:level0:g}) and 
(\ref{K+:f:level0:g}) 
replaced by 
$$
[K^+(z) , K^+(z')] = [K^{-(\al)}(z) , K^{-(\beta)}(z')] = 0,  
$$
\begin{align*}
& K^+(z) K^{-(\al)}(z_{\al})K^+(z)^{-1} K^{-(\al)}(z_{\al})^{-1}
\\ & = {{\theta(z-(z_{\al}) + \hbar h+\delta-\Delta)}
\over {\theta(z-(z_{\al}) - \hbar h+\delta-\Delta)}}
{{\theta(z-(z_{\al}) + \hbar h(K-1)+\delta-\Delta)}
\over {\theta(z-(z_{\al}) + \hbar h(K+1) +\delta-\Delta)}} , 
\end{align*}
$$
[e^{(\al)}(z_{\al}), f^{(\beta)}(w_{\beta}) ] = \delta_{\al\beta}
\left( \delta(z,(w_{\beta})) K^+((z_{\al})) - \delta(z,(w_{\beta});K) 
K^{-(\beta)}((w_{\beta}))^{-1} \right) ,  
$$
\begin{align*}
K^{-(\al)}(z_{\al}) e^{(\beta)}(w_{\beta}) K^{-(\al)}(z_{\al})^{-1} = 
{{\theta((z_{\al})-w_{\beta}-\hbar (K+1)h +\delta - \Delta)}\over
{\theta(z_{\al}-w_{\beta}-\hbar (K-1)h)}} e^{(\beta)}(w_{\beta}) , 
\end{align*}
where 
$$
\delta(z,w) = {{\pa_h\theta}\over\theta}(z-w+\delta-\Delta)
+
{{\pa_h\theta}\over\theta}(w-z+\delta-\Delta), 
$$
$$
\delta(z,w;K)=
{{\pa_h\theta}\over\theta}(z-w+\hbar Kh+\delta-\Delta) +
{{\pa_h\theta}\over\theta}(w-z+\hbar Kh+\delta-\Delta). 
$$ 

The construction of Prop. \ref{isom:i} and Cor. \ref{lancaster} 
of a deformation of the enveloping algebra of 
$\bar\G\otimes R_{S',\delta}$ in $U_{\hbar,h}\G_{S',\delta}$ 
may be extended to $U_{\hbar,h}\hat\G_{S',\delta}$.  

However, it seems difficult to construct a coproduct on 
$U_{\hbar,h}\hat\G_{S',\delta}$
because the maps $z\mapsto z_K$, where $z_K$ is the solution 
``close to $z$'' of
$\theta(z_K - z + \hbar K h + \delta - \Delta) = 0$, do not 
satisfy $(z_{K_1})_{K_2} = z_{K_1 + K_2}$.
\hfill \qed \medskip 
\end{remark}

\begin{remark} {\em Quantization of double extensions.}
A Hopf algebra $U_{\hbar}\wt\G_{S',\delta}$ quantizing the doubly 
extended quasi-Lie bialgebra 
$$
[(\bar \HH \otimes L_{S',\delta}^{\CC h}) \oplus 
(\bar \N_+ \otimes \cK_{S',\delta}) \oplus \CC D] 
\oplus
[(\bar \HH \otimes \wt\cO_{S',\delta}) \oplus 
(\bar \N_- \otimes \cK_{S',\delta}) \oplus \CC K] 
$$
may be obtained by replacing in the defining relations
of the above Rem., $q(z,w)$ by 
$$
q_{\pa}(z,w) = {{\theta(q^{\pa}z - w+\delta-\Delta)}\over
{\theta(z-q^{\pa}w+\delta-\Delta)}}
$$
and the shifts by $\hbar Kh$ by actions of $q^{K\pa}$ on 
the variables. One may then extend the construction of a 
subalgebra deforming the enveloping algebra of 
$\bar\G\otimes (\cO_{S'}\oplus\cK_{\delta}) \oplus\CC K$ 
to this situation. One may also construct an twist this Hopf 
structure obtain a quantization of
$$
[(\bar \HH \otimes R_{S',\delta}^{\CC h}) \oplus 
(\bar \N_+ \otimes \cK_{S',\delta}) \oplus \CC D] 
\oplus
[(\bar \HH \otimes \Sigma) \oplus 
(\bar \N_- \otimes \cK_{S',\delta}) \oplus \CC K] , 
$$
which will be isomorphic to the one defined in sect. \ref{general}. 
This is because the structure coefficient
$$
\wt q_{\delta}(z,w) = {{q_{\pa}(z,w)}\over{q_{\pa}(z,P)q_{\pa}(P,w)}}
$$
has the properties that it vanishes for $z = q^{\pa}w$, is 
single-valued for $z,w$ on $X$ and has its poles 
only for $z=w$ or $z,w$ in $\delta\cup \{P\}$.   

One may then construct in $U_{\hbar}\wt\G_{S',\delta}$, a deformation 
of the enveloping algebra of $R_{S',\delta}$ is the usual way.
Thm. \ref{isomorphisms} then implies that 
$U_{\hbar}\G'_{S',\delta}$ is isomorphic to the quotient 
$U_{\hbar,z^2dz}\G^{\prime \otimes (g-1)} \otimes DY(\SL_2)^{\prime
\otimes \card S'}$ by the identification of the central generators.  
\hfill \qed \medskip 
\end{remark}

\begin{remark}
If one replaces $\delta$ by an arbitrary effective divisor 
$\sum_i x_i$ in the definition of $U_{\hbar,h}\G_{S',\delta}$, 
the resulting algebra is no longer a flat deformation of $\wh
\SL_2$. For example, if one replaces the vertex relations by 
\begin{align*}
\theta(\sum_{i=1}^{g-1} x_{i} &  + (z_{\al})-(w_{\beta}) + 
\hbar h - \Delta) e^{(\al)}(z_{\al})e^{(\beta)}(w_{\beta})
\\ & = 
-\theta(\sum_{i=1}^{g-1} x_{i} + 
(w_{\beta}) - (z_{\al}) + \hbar h - \Delta)
e^{(\beta)}(w_{\beta})e^{(\al)}(z_{\al}), 
\end{align*}
one finds for $\al = \beta = x_i$, 
$$ 
(w_i(z_i-w_i) + \hbar \cdots) e^{(i)}(z_i)e^{(i)}(w_i) 
= (z_i(z_i-w_i) + \hbar\cdots) e^{(i)}(w_i)e^{(i)}(z_i), 
$$
which are relations for a flat deformation of the affinization 
of the Lie superalgebra ${\mathfrak osp}(2|1)$ (we owe this remark 
to B. Feigin). 
\end{remark}

\end{document}